\documentclass[11pt,leqno]{article}
\usepackage{amssymb}
\usepackage[mathscr]{eucal}
\usepackage{amsmath,amssymb,latexsym,theorem,bbm}
\usepackage{color,url}
\usepackage{appendix}

\newcommand{\CC}{\mathbb{C}}

\newcommand{\NN}{\mathbb{N}}
\newcommand{\QQ}{\mathbb{Q}}
\newcommand{\RR}{\mathbb{R}}

\newcommand{\ZZ}{\mathbb{Z}}

\newcommand{\bB}{{\boldsymbol{B}}}
\newcommand{\tb}{\widetilde{b}}
\newcommand{\tB}{\widetilde{B}}
\newcommand{\bc}{{\boldsymbol{c}}}

\newcommand{\bF}{{\boldsymbol{F}}}

\newcommand{\bM}{{\boldsymbol{M}}}

\newcommand{\bQ}{{\boldsymbol{Q}}}
\newcommand{\br}{{\boldsymbol{r}}}
\newcommand{\bR}{{\boldsymbol{R}}}

\newcommand{\bv}{{\boldsymbol{v}}}

\newcommand{\bx}{{\boldsymbol{x}}}

\newcommand{\bZ}{{\boldsymbol{Z}}}

\newcommand{\bmu}{{\boldsymbol{\mu}}}
\newcommand{\bbeta}{{\boldsymbol{\beta}}}
\newcommand{\bfeta}{{\boldsymbol{\eta}}}

\newcommand{\bpsi}{{\boldsymbol{\psi}}}
\newcommand{\btpsi}{{\widetilde{\bpsi}}}
\newcommand{\bzero}{{\boldsymbol{0}}}
\newcommand{\bone}{{\boldsymbol{1}}}

\newcommand{\cA}{{\mathcal A}}
\newcommand{\cB}{{\mathcal B}}

\newcommand{\cD}{{\mathcal D}}

\newcommand{\cI}{\mathcal{I}}

\newcommand{\cF}{{\mathcal F}}

\newcommand{\cN}{{\mathcal N}}
\newcommand{\cP}{{\mathcal P}}

\newcommand{\cV}{{\mathcal V}}

\newcommand{\cX}{{\mathcal X}}
\newcommand{\cY}{{\mathcal Y}}
\newcommand{\cZ}{{\mathcal Z}}
\newcommand{\cW}{{\mathcal W}}

\newcommand{\tcX}{\widetilde{\cX}}
\newcommand{\tcV}{\widetilde{\cV}}
\newcommand{\ttcV}{\widetilde{\tcV}}
\newcommand{\tcY}{\widetilde{\cY}}
\newcommand{\ttcY}{\widetilde{\tcY}}

\newcommand{\ttcX}{\widetilde{\tcX}}

\newcommand{\tcW}{\widetilde{\cW}}
\newcommand{\ttcW}{\widetilde{\tcW}}

\newcommand{\dd}{\mathrm{d}}
\newcommand{\ee}{\mathrm{e}}
\newcommand{\ii}{\mathrm{i}}

\newcommand{\cont}{\mathrm{cont}}

\DeclareMathOperator*{\argmax}{arg\,max}

\newcommand{\EE}{\operatorname{\mathbb{E}}}
\newcommand{\PP}{\operatorname{\mathbb{P}}}

\newcommand{\hb}{\widehat{b}}

\newcommand{\tbeta}{\widetilde{\beta}}

\newcommand{\tmu}{\widetilde{\mu}}

\newcommand{\tY}{\widetilde{Y}}

\newcommand{\tW}{\widetilde{W}}
\newcommand{\teta}{\widetilde{\eta}}
\newcommand{\hsigma}{\widehat{\sigma}}

\newcommand{\vare}{\varepsilon}

\renewcommand{\Re}{\operatorname{Re}}

\renewcommand{\mid}{\,|\,}

\renewcommand{\leq}{\leqslant}
\renewcommand{\geq}{\geqslant}

\newcommand{\stoch}{\stackrel{\PP}{\longrightarrow}}
\newcommand{\distr}{\stackrel{\cD}{\longrightarrow}}
\newcommand{\distre}{\stackrel{\cD}{=}}

\newcommand{\as}{\stackrel{{\mathrm{a.s.}}}{\longrightarrow}}
\newcommand{\ase}{\stackrel{{\mathrm{a.s.}}}{=}}

\newcommand{\bbone}{\mathbbm{1}}

\newcommand{\proofend}{\hfill\mbox{$\Box$}}

\numberwithin{equation}{section}

\theoremstyle{change} \theorembodyfont{\em}
\newtheorem{Lem}{Lemma.}[section]
\newtheorem{Thm}[Lem]{Theorem.}
\newtheorem{Pro}[Lem]{Proposition.}
\newtheorem{Cor}[Lem]{Corollary.}
\newtheorem{Def}[Lem]{Definition.}

\theorembodyfont{\rm}
\newtheorem{Rem}[Lem]{Remark.}
\newtheorem{Ex}[Lem]{Example.}

\begin{document}

\begin{center}
 {\bfseries\Large
  Asymptotic properties of maximum likelihood estimator \\[2mm]
   for the growth rate for a jump-type CIR process \\[2mm]
   based on continuous time observations} \\[7mm]
{\sc\large
 M\'aty\'as $\text{Barczy}^{*,\diamond}$, \ Mohamed $\text{Ben Alaya}^{**}$,}\\[2mm]
 {\sc\large Ahmed $\text{Kebaier}^{**}$ \ and \ Gyula $\text{Pap}^{***}$}
\end{center}

\vskip0.2cm

\noindent
 * MTA-SZTE Analysis and Stochastics Research Group,
   Bolyai Institute, University of Szeged,
   Aradi v\'ertan\'uk tere 1, H--6720 Szeged, Hungary.

\noindent
 ** Universit\'e Paris 13, Sorbonne Paris Cit\'e, LAGA, CNRS (UMR 7539),
    Villetaneuse, France.

\noindent
 *** Bolyai Institute, University of Szeged,
     Aradi v\'ertan\'uk tere 1, H--6720 Szeged, Hungary.

\noindent e--mails: barczy.matyas@inf.unideb.hu (M. Barczy), \\
\phantom{e--mails:\,} mba@math.univ-paris13.fr (M. Ben Alaya), \\
\phantom{e--mails:\,} kebaier@math.univ-paris13.fr (A. Kebaier), \\
\phantom{e--mails:\,} papgy@math.u-szeged.hu (G. Pap).

\noindent $\diamond$ Corresponding author.

\renewcommand{\thefootnote}{}
\vspace*{0.2cm}
\footnote{This research is supported by Laboratory of Excellence MME-DII, Grant no.\ ANR11-LBX-0023-01
 (\texttt{http://labex-mme-dii.u-cergy.fr/}).
M\'aty\'as Barczy was supported by the ''Magyar \'Allami E\"otv\"os
 \"Oszt\"ond\'{\i}j 2016'' Grant no.\ 75141 funded by the Tempus Public Foundation.
Ahmed Kebaier benefited from the support of the chair
 “Risques Financiers”, Fondation du Risque.}

\vspace*{-5mm}

\begin{abstract}
We consider a jump-type Cox--Ingersoll--Ross (CIR) process driven by a standard
 Wiener process and a subordinator, and we
 study asymptotic properties of the maximum likelihood estimator (MLE) for its growth rate.
We distinguish three cases: subcritical, critical and supercritical.
In the subcritical case we prove weak consistency and asymptotic normality,
 and, under an additional moment assumption, strong consistency as well.
In the supercritical case, we prove strong consistency and mixed normal (but non-normal)
 asymptotic behavior, while in the critical case, weak consistency and non-standard
 asymptotic behavior are described.
We specialize our results to so-called basic affine jump-diffusions as well.
Concerning the asymptotic behavior of the MLE in the supercritical case, we derive a
 stochastic representation of the limiting mixed normal distribution, where
 the almost sure limit of an appropriately scaled jump-type supercritical CIR process
 comes into play.
This is a new phenomenon, compared to the critical case, where a
 diffusion-type critical CIR process plays a role.
\end{abstract}

\section{Introduction}

Continuous state and continuous time branching processes with immigration, especially, the
 Cox--Ingersoll--Ross (CIR) process (introduced by Feller \cite{Fel} and Cox et
 al.\ \cite{CoxIngRos}) and its variants, play an important role in stochastics, and there is
 a wide range of applications of these processes in biology and financial mathematics as
 well.
In the framework of the famous Heston model, which is popular in finance, a CIR process can be interpreted as
 a stochastic volatility (or instantaneous variance) of a price process of an asset.
In this paper, we consider a jump-type CIR process driven by a standard
 Wiener process and a subordinator
 \begin{align}\label{jump_CIR}
  \dd Y_t = (a - b Y_t) \, \dd t + \sigma \sqrt{Y_t} \, \dd W_t + \dd J_t,
  \qquad  t \in [0, \infty) ,
 \end{align}
 with an almost surely non-negative initial value \ $Y_0$, \ where \ $a \in [0, \infty)$,
 \ $b \in \RR$, \ $\sigma \in (0, \infty)$, \ $(W_t)_{t \in [0, \infty)}$ \ is a
 1-dimensional standard Wiener process, and \ $(J_t)_{t\in[0,\infty)}$ \ is a subordinator
 (an increasing L\'evy process) with zero drift and with L\'evy measure \ $m$
 \ concentrating on \ $(0, \infty)$ \ such that
 \begin{align}\label{help_Levy}
  \int_0^\infty z \, m(\dd z) \in [0, \infty) ,
 \end{align}
 that is,
 \begin{equation}\label{LK}
  \EE(\ee^{u J_t})
  = \exp\left\{ t \int_0^\infty (\ee^{uz} - 1) \, m(\dd z) \right\}
 \end{equation}
 for any \ $t \in [0, \infty)$ \ and for any complex number \ $u$ \ with
 \ $\Re(u) \in (-\infty, 0]$, \ see, e.g., Sato \cite[proof of Theorem 24.11]{Sat}.
We suppose that \ $Y_0$, \ $(W_t)_{t \in [0, \infty)}$ \ and \ $(J_t)_{t\in[0,\infty)}$
 \ are independent.
Note that the moment condition \eqref{help_Levy} implies that \ $m$ \ is a L\'evy measure
 (since \ $\min(1, z^2) \leq z$ \ for \ $z \in (0, \infty)$).
Moreover, the subordinator \ $J$ \ has sample paths of bounded variation on every compact time interval almost surely,
 see, e.g., Sato \cite[Theorem 21.9]{Sat}.
We point out that the assumptions assure that there is a (pathwise) unique strong solution
 of the SDE \eqref{jump_CIR} with
 \ $\PP(\text{$Y_t \in [0, \infty)$ for all $t \in [0, \infty)$}) = 1$
 \ (see Proposition \ref{Pro_jump_CIR}).
In fact, \ $(Y_t)_{t\in[0,\infty)}$ \ is a special continuous state and continuous time
 branching process with immigration (CBI process), see Proposition \ref{Pro_jump_CIR}.

In the present paper, we focus on parameter estimation for the jump-type CIR process \eqref{jump_CIR}
 in critical and supercritical cases \ ($b=0$ \ and \ $b\in(-\infty,0)$, \ respectively),
 which have not been addressed in previous research.
We also study the subcritical case \ $(b\in(0,\infty))$ \ and we get results extending those of Mai
 \cite[Theorem 4.3.1]{Mai} in several aspects: we do not suppose the ergodicity of
 the process \ $Y$ \ and we make explicit the expectation of the unique stationary
 distribution of \ $Y$ \ in the limit law in Theorem \ref{Thm_MLEb_subcritical}.
However, we note that some points in Mai's approach \cite[Sections 3.3 and 4.3]{Mai} should be
 corrected concerning the expressions of the likelihood ratio (Mai \cite[formula (3.10)]{Mai})
 and the maximum likelihood estimator (MLE) of \ $b\in\RR$ \ (Mai \cite[formula (4.23)]{Mai}),
 see our results in Propositions \ref{RNb} and \ref{LEMMA_MLEb_exist}, respectively.
Supposing that \ $a \in [0, \infty)$, \ $\sigma \in (0, \infty)$ \ and the measure
 \ $m$ \ are known, we study the asymptotic properties of the MLE of \ $b \in \RR$ \
 based on continuous time observations \ $(Y_t)_{t\in[0,T]}$ \ with \ $T \in (0, \infty)$,
 \ starting the process \ $Y$ \ from some known non-random initial value \ $y_0 \in [0, \infty)$.
\ It will turn out that for the calculation of the MLE of \ $b$, \ one does not need to know
 the value of the parameter \ $\sigma$ \ and the measure \ $m$, \ see \eqref{MLEb}.
We have restricted ourselves to studying the MLE of \ $b$ \ supposing that \ $a$ \ is known,
 since in order to describe the asymptotic behavior of the MLE of \ $a$ \  supposing
 that \ $b$ \ is known or the joint MLE of \ $(a, b)$, \ one has to find, for
 instance, the limiting behavior of \ $\int_0^t \frac{1}{Y_s} \, \dd s$ \ as
 \ $t \to \infty$, \ which seems to be a hard task even in the subcritical case.
In general, we would need an explicit formula for the Laplace transform
 of \ $\int_0^t \frac{1}{Y_s} \, \dd s$, \ $t\in\RR_+$, \ which is not known up to our knowledge.
This can be a topic of further research.

Studying asymptotic properties of various kinds of estimators for the drift parameters of the
 CIR process and its variants has a long history, but most of the existing results refer to
 the original (diffusion-type) CIR process.
Overbeck \cite{Ove} studied the MLE of the drift parameters of the original CIR process based
 on continuous time observations, and later on, Ben Alaya and Kebaier \cite{BenKeb1},
 \cite{BenKeb2} completed the results of Overbeck \cite{Ove} giving explicit forms of the
 joint Laplace transforms of the building blocks of the MLE in question as well.
Another type of estimator, so-called conditional least squares estimator (LSE) has also
 been investigated for the drift parameters for the original CIR process, see, e.g.,
 Overbeck and Ryd\'en \cite{OveRyd}.
For a generalization of the original CIR process, namely, for a CIR model driven by a stable
 noise instead of a standard Wiener process (also called a stable CIR model) Li and Ma
 \cite{LiMa} described the asymptotic
 behaviour of the (weighted) conditional LSE of the drift parameters of this model based on a
 discretely observed low frequency data set in the subcritical case.
For a CBI process, being a generalization of a (stable) CIR process, Huang et al.
 \cite{HuaMaZhu} studied the asymptotics of the weighted conditional LSE of the drift
 parameters of the model based on low frequency discrete time observations under second order
 moment assumptions of the branching and immigration mechanisms of the CBI process in
 question.

Note that we have \ $\EE(J_t) = t \int_0^\infty z \, m(\dd z) \in [0, \infty)$, \ and
 \ $\EE(Y_t) \in [0, \infty)$ \ (see Proposition \ref{Pro_moments}).
Moreover, \ $(J_t)_{t\in[0,\infty)}$ \ is a compound Poisson process if and only if
 \ $m(\RR) = m((0,\infty)) \in [0, \infty)$, \ see, e.g., Sato \cite[Examples 8.5]{Sat}.
If \ $m((0,\infty)) = \infty$, \ then its jump intensity is infinity, i.e., almost surely,
 the jump times are infinitely many, countable and dense in \ $[0, \infty)$, \ and if
 \ $m((0,\infty)) \in (0, \infty)$, \ then, almost surely,
 there are finitely many jump times on every compact intervals yielding that the jump times are infinitely
 many and countable in increasing order, and the jump intensity is \ $m((0,\infty))$, \ i.e.,
 the first jump time has an exponential distribution with mean \ $1/m((0,\infty))$, \ and
 the distribution of the jump size is \ $m(\dd z)/m((0,\infty))$
 \ (see, e.g., Sato \cite[Theorem 21.3]{Sat}).
The case \ $m((0,\infty)) = 0$ \ corresponds to the usual CIR process.
Our forthcoming results will cover both cases \ $m((0,\infty))\in [0,\infty)$ \ and \ $m((0,\infty)) = \infty$.

In case of \ $b \in (0, \infty)$, \ $Y_t$ \ converges in law as \ $t \to \infty$ \ to its
 unique stationary distribution \ $\pi$ \ (see Theorem \ref{Ergodicity}).
This follows from a general result for CBI processes which has been announced without proof
 in Pinsky \cite{Pin} and a proof has been given in Li
 \cite[Theorem 3.20 and the paragraph after Corollary 3.21]{Li}, see also Keller-Ressel and
 Steiner \cite{KRSte}, Keller-Ressel \cite{KR}, and Keller-Ressel and Mijatovi\'c
 \cite[Theorem 2.6]{KRMij}.
The mean \ $\int_0^\infty y \, \pi(\dd y) \in [0, \infty)$ \ of the unique stationary
 distribution is the so-called long variance (long run average price variance, i.e., the
 limit of \ $\EE(Y_t)$ \ as \ $t \to \infty$, \ see \eqref{expectation_sub} and
 \eqref{help_stac_expectation}), \ $b$ \ is the rate at which \ $\EE(Y_t)$ \ reverts to
 \ $\int_0^\infty y \, \pi(\dd y)$ \ as \ $t \to \infty$ \ (speed of adjustment, since
 \ $\EE(Y_t)
    = \int_0^\infty y \, \pi(\dd y)
      + \ee^{-bt} \bigl(\EE(Y_0) - \int_0^\infty y \, \pi(\dd y)\bigr)$
 \ for all \ $t \in [0, \infty)$, \ see \eqref{expectation} and
 \eqref{help_stac_expectation}).
Under \ $a \in (0,\infty)$, \ the moment condition \eqref{help_Levy} and the extra moment
 condition
 \begin{equation}\label{EXTRA}
  \int_0^1 z \log\left(\frac{1}{z}\right) m(\dd z) < \infty ,
 \end{equation}
 Jin et al.\ \cite{JinRudTra2} established an explicit positive lower bound of the transition
 densities of \ $(Y_t)_{t\in[0,\infty)}$, \ and based on this result, they showed the
 existence of a Foster--Lyapunov function and derived exponential ergodicity for
 \ $(Y_t)_{t\in[0,\infty)}$ \ (see Theorem \ref{Ergodicity}).
Comparing the moment conditions \eqref{help_Levy} and \eqref{EXTRA}, note that the
 integrability of \ $z \log\left(\frac{1}{z}\right)$ \ on the interval \ $(0, \ee^{-1})$
 \ yields that of \ $z$ \ on the same interval.

In case of \ $b = 0$, \ if \ $a + \int_0^\infty z \, m(\dd z) \in (0,\infty)$, \ then
 \ $\lim_{t \to \infty}\EE(Y_t) = \infty$ \ such that
 \ $\lim_{t\to\infty} t^{-1} \EE(Y_t) \in (0, \infty)$, \ and,
 in case of \ $b\in(-\infty,0)$, \ if \ $\EE(Y_0)\in(0,\infty)$ \ or
 \ $a + \int_0^\infty z \, m(\dd z) \in (0,\infty)$ \ (which rule out the case that \ $Y$ \ is identically zero),
 then \ $\lim_{t \to \infty}\EE(Y_t) = \infty$ \ such that \ $\lim_{t\to\infty} \ee^{bt} \EE(Y_t) \in (0, \infty)$,
 \ hence the parameter \ $b$ \ can always be interpreted as the growth rate, see Proposition
 \ref{Pro_moments}.

The jump-type CIR process in \eqref{jump_CIR} includes the so-called basic affine
 jump-diffusion (BAJD) as a special case, in which the drift takes the form
 \ $\kappa (\theta - Y_t)$ \ with some \ $\kappa \in (0, \infty)$ \ and
 \ $\theta \in [0, \infty)$, \ and the L\'evy process \ $(J_t)_{t\in[0,\infty)}$ \ is a
 compound Poisson process with exponentially distributed jump sizes, namely,
 \begin{align}\label{Levy_cond}
  m(\dd z) =  c \lambda \ee^{-\lambda z} \bbone_{(0,\infty)}(z) \, \dd z
 \end{align}
 with some constants \ $c \in [0, \infty)$ \ and \ $\lambda \in (0, \infty)$.
\ Note that the measure \ $m$ \ given by \eqref{Levy_cond} satisfies \eqref{help_Levy} and
 \eqref{EXTRA}, and, for the compound Poisson process in question, the first jump time has an exponential
 distribution with parameter \ $c$ \ and the distribution of the jump size is exponential with parameter \ $\lambda$.
\ Indeed, in this special case
 \begin{align*}
  \int_0^\infty z \, m(\dd z) = c \int_0^\infty z \lambda \ee^{-\lambda z} \, \dd z
                              = \frac{c}{\lambda} \in [0, \infty) ,
 \end{align*}
 and
 \begin{align*}
  \int_0^1 z\log\left(\frac{1}{z}\right) m(\dd z)
  \leq c \lambda \int_0^1 z\log\left(\frac{1}{z}\right) \, \dd z
  = c \lambda \int_0^\infty u \ee^{-2u} \, \dd u
  = \frac{c\lambda}{4} \in [0, \infty).
 \end{align*}
For describing the dynamics of default intensity, the BAJD was introduced by Duffie and
 G\^arleanu \cite{DufGar}.
Filipovi\'c \cite{Fil} and Keller-Ressel and Steiner \cite{KRSte} used the BAJD as a
 short-rate model.

The paper is organized as follows.
In Section \ref{Prel}, we prove that the SDE \eqref{jump_CIR} has a pathwise unique strong
 solution (under some appropriate conditions), see Proposition \ref{Pro_jump_CIR}.
We describe the asymptotic behaviour of the first moment of \ $(Y_t)_{t\in[0,\infty)}$,
 \ and, based on it, we introduce a classification of jump-type CIR processes given by the
 SDE \eqref{jump_CIR}, see Proposition \ref{Pro_moments} and Definition
 \ref{Def_criticality}.
Namely, we call \ $(Y_t)_{t\in[0,\infty)}$ \ subcritical, critical or supercritical if
 \ $b \in (0,\infty)$, \ $b = 0$, \ or \ $b \in (-\infty, 0)$, \ respectively.
We recall a result about the existence of a unique stationary distribution and exponential
 ergodicity for the process \ $(Y_t)_{t\in[0,\infty)}$ \ given by the equation
 \eqref{jump_CIR}, see Theorem \ref{Ergodicity}.
In Remark \ref{Rem_Grigelionis}, we derive a Grigelionis representation for the process
 \ $(Y_t)_{t\in[0,\infty)}$.
\ Further, we explain why we do not estimate the parameter \ $\sigma$, \ see Remark
 \ref{Thm_MLE_cons_sigma}.
Next we drive explicit formulas for the Laplace transform of
 \ $\bigl(Y_t, \int_0^t Y_s \, \dd s\bigr)$ \ in Section \ref{section_Laplace}, together
 with some examples for the BAJD process.
Here we use the fact that \ $\bigl(Y_t, \int_0^t Y_s \, \dd s\bigr)_{t\in[0,\infty)}$ \ is a 2-dimensional CBI process
 following also from Keller-Ressel \cite[Theorem 4.10]{KR2} or Filipovi\'{c} et al. \cite[paragraph before Theorem 4.3]{FilMaySch}.
For completeness, we note that Keller-Ressel \cite[Theorem 4.10]{KR2} derived a
 formula for the joint Laplace transform of a regular affine process and its integrated process
 containing the solutions of Riccati-type differential equations,
 and Jiao et al. \cite[Proposition 4.3]{JiaMaSco} derived a formula for that of a general CBI process and its integrated process.
We point out that our proof of technique of Theorem \ref{Thm_Laplace_joint}
 is different from those of Keller-Ressel \cite[Theorem 4.10]{KR2} and Jiao et al. \cite[Proposition 4.3]{JiaMaSco},
 and we make the solutions of Riccati-type differential equations explicit in case of \ $(Y_t)_{t\in[0,\infty)}$.
\ Section \ref{section_EUMLE} is devoted to study the existence and uniqueness of the
 MLE \ $\hb_T$ \ of \ $b$ \ based on observations \ $(Y_t)_{t\in[0,T]}$ \ with
 \ $T \in (0, \infty)$.
\ We derive an explicit formula for \ $\hb_T$ \ as well, see \eqref{MLEb}.
Sections \ref{section_MLE_subcritical}, \ref{section_MLE_critical} and
 \ref{section_MLE_supercritical} are devoted to study asymptotic behaviour of the MLE of
 \ $b$ \ for subcritical, critical and supercritical jump-type CIR models, respectively.
In Section \ref{section_MLE_subcritical}, we show that in the subcritical case, the MLE of
 \ $b$ \ is asymptotically normal with the usual square root scaling \ $T^{1/2}$ \ (especially,
 it is weakly consistent), but unfortunately, the asymptotic variance depends on the unknown parameters
 \ $a$ \ and \ $m$, \ as well.
To get around this problem, we also replace the deterministic scaling \ $T^{1/2}$ \ by the
 random scaling \ $\frac{1}{\sigma} \bigl(\int_0^T Y_s\,\dd s\bigr)^{1/2}$ \ with the
 advantage that the MLE of \ $b$ \ with this scaling is asymptotically standard normal.
Under the extra moment condition \eqref{EXTRA}, we prove strong consistency as well.
In Section \ref{section_MLE_critical}, we describe the (non-normal) asymptotic behaviour of
 the MLE of \ $b$ \ in the critical case both with the deterministic scaling \ $T$ \ and with
 the random scaling \ $\frac{1}{\sigma} \bigl(\int_0^T Y_s\,\dd s\bigr)^{1/2}$.
In Section \ref{section_MLE_supercritical}, for the supercritical case, we prove that the
 MLE of \ $b$ \ is strongly consistent, and it is asymptotically mixed normal with the
 deterministic scaling \ $\ee^{-bT/2}$, \ and it is asymptotically standard normal with the
 random scaling \ $\frac{1}{\sigma} \bigl(\int_0^T Y_s\,\dd s\bigr)^{1/2}$.
\ We close the paper with Appendices, where we prove a comparison theorem for the SDE
 \eqref{jump_CIR} in the jump process \ $J$, \ we recall certain sufficient conditions for
 the absolute continuity of probability measures induced by semimartingales together with a
 representation of the Radon--Nikodym derivative (Appendix \ref{App_LR}) and some limit
 theorems for continuous local martingales (Appendix \ref{App_clm}).

Finally, we summarize the novelties of the paper.
We point out that only few results are available for parameter estimation for jump-type CIR processes,
 see Mai \cite[Section 4.3]{Mai} (MLE for subcritical case), Huang, Ma and Zhu \cite{HuaMaZhu} and
 Li and Ma \cite{LiMa} (conditional LSE).
Concerning the asymptotic behavior of the MLE in the subcritical case, we use an explicit
 formula for the Laplace transform of \ $\int_0^t Y_s \, \dd s$ \ to derive stochastic
 convergence of \ $\frac{1}{t} \int_0^t Y_s \, \dd s$ \ as \ $t \to \infty$, \ and we prove
 asymptotic normality avoiding ergodicity, see Theorem \ref{Thm_MLEb_subcritical}.
Further, in the supercritical case, we derive a stochastic representation in Theorem
 \ref{Thm_supercritical_convergence} of the limiting mixed normal distribution given in
 Theorem \ref{Thm_MLE_supercritical}, where the almost sure limit of an appropriately scaled
 jump-type supercritical CIR process comes into play.
This is a new phenomenon, compared to the critical case in Theorem \ref{Thm_MLE_critical_spec},
 where a diffusion-type critical CIR process plays a role.
We remark that for all \ $b \in \RR$,
 \ $\frac{1}{\sigma} \bigl(\int_0^T Y_s\,\dd s\bigr)^{1/2} (\hb_T - b)$ \ converges in
 distribution as \ $T \to \infty$, \ and the limit distribution is standard normal for the
 non-critical cases, while it is non-normal for the critical case (given explicitly in
 Theorem \ref{Thm_MLE_critical_spec}).
Hence we have a kind of unified theory.

\section{Preliminaries}
\label{Prel}

Let \ $\NN$, \ $\ZZ_+$, \ $\RR$, \ $\RR_+$, \ $\RR_{++}$, \ $\RR_-$, \ $\RR_{--}$ \ and
 \ $\CC$ \ denote the sets of positive integers, non-negative integers,
 real numbers, non-negative real numbers, positive real numbers, non-positive
 real numbers, negative real numbers and complex numbers, respectively.
For \ $x , y \in \RR$, \ we will use the notations \ $x \land y := \min(x, y)$
 \ and \ $x \lor y := \max(x, y)$.
\ The integer part of a real number \ $x \in \RR$ \ is denoted by \ $\lfloor x \rfloor$.
\ By \ $\|x\|$ \ and \ $\|A\|$, \ we denote the Euclidean norm of a vector
 \ $x \in \RR^d$ \ and the induced matrix norm of a matrix
 \ $A \in \RR^{d \times d}$, \ respectively.
By \ $\cB(\RR_+)$, \ we denote the Borel $\sigma$-algebra on \ $\RR_+$.
\ We will denote the convergence in probability, in distribution and almost surely, and
 equality in distribution and almost surely by \ $\stoch$, \ $\distr$, \ $\as$,
 \ $\distre$ \ and \ $\ase$, \ respectively.

Let \ $\bigl(\Omega, \cF, \PP\bigr)$ \ be a probability space.
By \ $C^2_c(\RR_+, \RR)$ \ and \ $C^{\infty}_c(\RR_+, \RR)$,
 \ we denote the set of twice continuously differentiable real-valued
 functions on \ $\RR_+$ \ with compact support, and the set of
 infinitely differentiable real-valued functions on \ $\RR_+$ \ with
 compact support, respectively.

The next proposition is about the existence and uniqueness of a strong
 solution of the SDE \eqref{jump_CIR} stating also that \ $Y$ \ is a CBI process.

\begin{Pro}\label{Pro_jump_CIR}
Let \ $\eta_0$ \ be a random variable independent of \ $(W_t)_{t\in\RR_+}$ \ and
 \ $(J_t)_{t\in\RR_+}$ \ satisfying \ $\PP(\eta_0 \in \RR_+) = 1$ \ and
 \ $\EE(\eta_0) < \infty$.
\ Then for all \ $a \in \RR_+$, \ $b \in \RR$, \  $\sigma \in \RR_{++}$ \ and L\'evy
 measure \ $m$ \ on \ $\RR_{++}$ \ satisfying \eqref{help_Levy}, there is a pathwise
 unique strong solution \ $(Y_t)_{t\in\RR_+}$ \ of the SDE \eqref{jump_CIR} such that
 \ $\PP(Y_0 = \eta_0) = 1$ \ and
 \ $\PP(\text{$Y_t \in \RR_+$ \ for all \ $t \in \RR_+$}) = 1$.
\ Moreover, \ $(Y_t)_{t\in\RR_+}$ \ is a CBI process having branching mechanism
 \[
   R(u)= \frac{\sigma^2}{2} u^2 - bu , \qquad
   \text{$u \in \CC$ \ with \ $\Re(u) \leq 0$,}
 \]
 and immigration mechanism
 \[
   F(u) = a u + \int_0^\infty (\ee^{uz} - 1) \, m(\dd z) ,
   \qquad \text{$u \in \CC$ \ with \ $\Re(u) \leq 0$.}
 \]
Further, the infinitesimal generator of \ $Y$ \ takes the form
 \begin{align}\label{infgen}
  \begin{split}
   (\cA f)(y) = (a - by) f'(y) + \frac{1}{2} y \sigma^2 f''(y)
                + \int_0^\infty (f(y+z) - f(y)) \, m(\dd z) ,
  \end{split}
 \end{align}
 where \ $y \in \RR_+$, \ $f \in C^2_c(\RR_+, \RR)$, \ and \ $f'$ \ and \ $f''$ \ denote
 the first and second order partial derivatives of \ $f$.

If, in addition, \ $\PP(\eta_0 \in \RR_{++}) = 1$ \ or \ $a \in \RR_{++}$, \ then
 \ $\PP\bigl(\int_0^t Y_s \, \dd s \in \RR_{++}\bigr) = 1$ \ for all \ $t \in \RR_{++}$.%
\end{Pro}

\noindent{\bf Proof.}
The L\'evy--It\^{o}'s representation of \ $J$ \ takes the form
 \begin{align}\label{J_Levy_Ito}
  J_t = \int_0^t \int_0^\infty z \, \mu^J(\dd s,\dd z), \qquad t \in \RR_+ ,
 \end{align}
 where
 \ $\mu^J(\dd s,\dd z)
    := \sum_{u\in\RR_+}
        \bbone_{\{\Delta J_u\ne 0\}} \, \vare_{(u,\Delta J_u)}(\dd s,\dd z)$
 \ is the integer-valued Poisson random measure on \ $\RR_{++}^2$ \ associated with the
 jumps \ $\Delta J_u := J_u - J_{u-}$, \ $u \in \RR_{++}$, \ $\Delta J_0 := 0$, \ of the
 process \ $J$, \ and \ $\vare_{(u,x)}$ \ denotes the Dirac measure at the point
 \ $(u, x) \in \RR_+^2$, \ see, e.g., Sato \cite[Theorem 19.2]{Sat}.
Consequently, the SDE \eqref{jump_CIR} can be rewritten in the form
 \begin{align}\label{jump_CIR_rewritten}
  \begin{split}
   Y_t &= Y_0 + \int_0^t (a - bY_s) \, \dd s
          + \int_0^t \sigma \sqrt{Y_s} \, \dd W_s + J_t\\
       &= Y_0 + \int_0^t (a - bY_s) \, \dd s + \int_0^t \sigma \sqrt{Y_s} \, \dd W_s
          + \int_0^t \int_0^\infty z \, \mu^J(\dd s,\dd z) ,
   \qquad t \in \RR_+.
  \end{split}
 \end{align}
Equation \eqref{jump_CIR_rewritten} is a special case of the equation (6.6) in Dawson
 and Li \cite{DawLi}, and Theorem 6.2 in Dawson and Li \cite{DawLi} implies that for any
 initial value \ $\eta_0$ \ with \ $\PP(\eta_0 \in \RR_+) = 1$ \ and \ $\EE(\eta_0)<\infty$,
 \ there exists a pathwise unique non-negative strong solution satisfying
 \ $\PP(Y_0 = \eta_0) = 1$ \ and
 \ $\PP(\text{$Y_t \in \RR_+$ \ for all \ $t \in \RR_+$}) = 1$.
\ Let \ $(Y'_t)_{t\in\RR_+}$ \ be a pathwise unique non-negative strong solution of the SDE
 \[
   \dd Y'_t = (a - b Y'_t) \, \dd t + \sigma \sqrt{Y'_t} \, \dd W_t ,
   \qquad t \in \RR_+ ,
 \]
 such that \ $\PP(Y'_0 = \eta_0) = 1$.
\ Applying the comparison Theorem \ref{Pro_comparison_jump_CIR}, we obtain
 \begin{align}\label{help_int_finite}
   \PP(\text{$Y_t \geq Y_t'$ \ for all \ $t \in \RR_+$}) = 1 .
 \end{align}
Further, if \ $\PP(\eta_0 \in \RR_{++}) = 1$ \ or \ $a \in \RR_{++}$, \ then
 \ $\PP\bigl(\int_0^t Y_s' \, \dd s \in \RR_{++}\bigr) = 1$ \ for all \ $t \in \RR_{++}$.
\ Indeed, if \ $\omega \in \Omega$ \ is such that \ $[0, t] \ni u \mapsto Y_u'(\omega)$ \ is
 continuous and \ $Y_v'(\omega) \in \RR_+$ \ for all \ $v \in\RR_+$, \ then we have
 \ $\int_0^t Y_s'(\omega) \, \dd s = 0$ \ if and only if \ $Y_s'(\omega) = 0$ \ for all
 \ $s \in [0, t]$.
\ Using the method of the proof of Theorem 3.1 in Barczy et.\ al \cite{BarDorLiPap}, we get
 \ $\PP\bigl(\int_0^t Y_s' \, \dd s = 0\bigr) = 0$, \ $t \in \RR_+$, \ as desired.
Since \ $(Y_s)_{s\in[0,t]}$ \ has c\`adl\`ag, hence bounded sample paths almost surely (see,
 e.g., Billingsley \cite[(12.5)]{Bil}), using \eqref{help_int_finite}, we conclude
 \ $\PP\bigl(\int_0^t Y_s \, \dd s \in \RR_{++}\bigr) = 1$ \ for all \ $t \in \RR_{++}$.

The form of the infinitesimal generator \eqref{infgen} readily follows by (6.5) in Dawson
 and Li \cite{DawLi}.
Further, Theorem 6.2 in Dawson and Li \cite{DawLi} also implies that \ $Y$ \ is a continuous
 state and continuous time branching process with immigration having branching and
 immigration mechanisms given in the Proposition.
\proofend

Next we present a result about the first moment of \ $(Y_t)_{t\in\RR_+}$.

\begin{Pro}\label{Pro_moments}
Let \ $a \in \RR_+$, \ $b \in \RR$, \ $\sigma \in \RR_{++}$, \ and let \ $m$ \ be a L\'evy
 measure on \ $\RR_{++}$ \ satisfying \eqref{help_Levy}.
Let \ $(Y_t)_{t\in\RR_+}$ \ be the unique strong solution of the SDE
 \eqref{jump_CIR} satisfying \ $\PP(Y_0 \in \RR_+) = 1$ \ and \ $\EE(Y_0) < \infty$.
\ Then
 \begin{align}\label{expectation}
  \EE(Y_t)
  = \begin{cases}
     \ee^{-bt} \EE(Y_0)
     + \bigl(a + \int_0^\infty z \, m(\dd z)\bigr) \frac{1-\ee^{-bt}}{b}
      & \text{if \ $b \ne 0$,} \\[1mm]
     \EE(Y_0) + \bigl(a + \int_0^\infty z \, m(\dd z)\bigr) t
      & \text{if \ $b = 0$,}
    \end{cases}
  \qquad t \in \RR_+ .
 \end{align}
Consequently, if \ $b \in \RR_{++}$, \ then
 \begin{equation}\label{expectation_sub}
  \lim_{t\to\infty} \EE(Y_t)
  = \biggl(a + \int_0^\infty z \, m(\dd z)\biggr) \frac{1}{b} ,
 \end{equation}
 if \ $b = 0$, \ then
 \[
   \lim_{t\to\infty} t^{-1} \EE(Y_t) = a + \int_0^\infty z \, m(\dd z) ,
 \]
 if \ $b \in \RR_{--}$, \ then
 \[
   \lim_{t\to\infty} \ee^{bt} \EE(Y_t)
   = \EE(Y_0) - \biggl(a + \int_0^\infty z \, m(\dd z) \biggr) \frac{1}{b} .
 \]
\end{Pro}

\noindent{\bf Proof.}
By Proposition \ref{Pro_jump_CIR}, \ $(Y_t)_{t\in\RR_+}$ \ is CBI process with an
 infinitesimal generator given in \eqref{infgen}.
By the notations of Barczy et al.\ \cite{BarLiPap2}, this CBI process has parameters
 \ $(d, c, \beta, B, \nu, \mu)$, \ where \ $d = 1$, \ $c = \frac{1}{2} \sigma^2$,
 \ $\beta = a$, \ $B = -b$, \ $\nu = m$ \ and \ $\mu = 0$.
\ Since \ $\EE(Y_0) < \infty$ \ and the moment condition
 \ $\int_{\RR\setminus\{0\}} |z| \bbone_{\{|z|\geq1\}} \, \nu(\dd z)
  = \int_1^\infty z \, m(\dd z) < \infty$ \ holds (due to \eqref{help_Levy}), we
 may apply formula (3.1.11) in Li \cite{Li2}
 or Lemma 3.4 in Barczy et al.\ \cite{BarLiPap2} with the choices \ $\tB = B = -b$ \ and
 \[
   \tbeta
   = \beta + \int_{\RR\setminus\{0\}} z \, \nu(\dd z)
   = \beta + \int_0^\infty z \, m(\dd z)
   \in\RR_+,
 \]
 yielding that
 \[
   \EE(Y_t) = \ee^{t\tB} \EE(Y_0) + \left(\int_0^t \ee^{u\tB} \, \dd u\right) \tbeta .
 \]
This implies \eqref{expectation} and the other parts of the assertion.
\proofend

Based on the asymptotic behavior of the expectations \ $\EE(Y_t)$ \ as \ $t \to \infty$,
 \ we introduce a classification of jump-type CIR model driven by a subordinator given by
 the SDE \eqref{jump_CIR}.

\begin{Def}\label{Def_criticality}
Let \ $a \in \RR_+$, \ $b \in \RR$, \ $\sigma \in \RR_{++}$, \ and let \ $m$ \ be a L\'evy
 measure on \ $\RR_{++}$ \ satisfying \eqref{help_Levy}.
Let \ $(Y_t)_{t\in\RR_+}$ \ be the unique strong solution of the SDE \eqref{jump_CIR}
 satisfying \ $\PP(Y_0 \in \RR_+) = 1$ \ and \ $\EE(Y_0) < \infty$.
\ We call \ $(Y_t)_{t\in\RR_+}$ \ subcritical, critical or supercritical if
 \ $b \in \RR_{++}$, \ $b = 0$ \ or \ $b \in \RR_{--}$, \ respectively.
\end{Def}

In the subcritical case, the following result states the existence of a unique
 stationary distribution and the exponential ergodicity for the process
 \ $(Y_t)_{t\in\RR_+}$, \ see Pinsky \cite{Pin},
 Li \cite[Theorem 3.20 and the paragraph after Corollary 3.21]{Li},
 Keller-Ressel and Steiner \cite{KRSte}, Keller-Ressel \cite{KR}, Keller-Ressel and
 Mijatovi\'c \cite[Theorem 2.6]{KRMij} and Jin et al. \cite[Theorem 1]{JinRudTra2}.
As a consequence, according to the discussion after Proposition 2.5 in Bhattacharya
 \cite{Bha}, one also obtains a strong law of large numbers for \ $(Y_t)_{t\in\RR_+}$.

\begin{Thm}\label{Ergodicity}
Let \ $a \in \RR_+$, \ $b \in \RR_{++}$, \ $\sigma \in \RR_{++}$, \ and let \ $m$ \ be
 a L\'evy measure on \ $\RR_{++}$ \ satisfying \eqref{help_Levy}.
Let \ $(Y_t)_{t\in\RR_+}$ \ be the unique strong solution of the SDE
 \eqref{jump_CIR} satisfying \ $\PP(Y_0 \in \RR_+) = 1$ \ and \ $\EE(Y_0) < \infty$.%
 \renewcommand{\labelenumi}{{\rm(\roman{enumi})}}
 \begin{enumerate}
  \item
   Then \ $(Y_t)_{t\in\RR_+}$ \ converges in law to its unique stationary distribution
   \ $\pi$ \ given by
    \[
      \int_0^\infty \ee^{uy} \, \pi(\dd y)
      = \exp\biggr\{\int_u^0 \frac{F(v)}{R(v)} \, \dd v\biggl\}
      = \exp\biggr\{\int_u^0
         \frac{av+\int_0^\infty(\ee^{vz}-1)\,m(\dd z)}
              {\frac{\sigma^2}{2}\,v^2-bv} \, \dd v\biggl\}
    \]
    for \ $u \in \RR_-$.
   \ Moreover, \ $\pi$ \ has a finite expectation given by
    \begin{align}\label{help_stac_expectation}
     \int_0^\infty y \, \pi(\dd y)
     = \biggl(a + \int_0^\infty z \, m(\dd z)\biggr) \frac{1}{b} \in \RR_+ .
    \end{align}
 \item
  If, in addition, \ $a \in \RR_{++}$ \ and the extra moment condition \eqref{EXTRA} holds,
   then the process \ $(Y_t)_{t\in\RR_+}$ \ is \ exponentially ergodic, namely, there exist
   constants \ $\beta \in (0, 1)$ \ and \ $C \in \RR_{++}$ \ such that
   \[
     \|\PP_{Y_t \mid Y_0 = y} - \pi\|_{\mathrm{TV}}
     \leq C (y + 1) \beta^t , \qquad t \in \RR_+ , \qquad y \in \RR_+ ,
   \]
   where \ $\|\mu\|_{\mathrm{TV}}$ \ denotes the total-variation norm of a signed measure
   \ $\mu$ \ on \ $\RR_+$ \ defined by
   \ $\|\mu\|_{\mathrm{TV}} := \sup_{A\in\cB(\RR_+)} |\mu(A)|$, \ and
   \ $\PP_{Y_t \mid Y_0 = y}$ \ is the conditional distribution of \ $Y_t$ \ with respect to
   the condition \ $Y_0 = y$.
   \ Moreover, for all Borel measurable functions \ $f : \RR_+ \to \RR$ \ with
   \ $\int_0^\infty |f(y)| \, \pi(\dd y) < \infty$, \ we have
   \begin{equation}\label{ergodic}
    \frac{1}{T} \int_0^T f(Y_s) \, \dd s \as \int_0^\infty f(y) \, \pi(\dd y) \qquad
    \text{as \ $T \to \infty$.}
   \end{equation}
 \end{enumerate}
\end{Thm}

\begin{Rem}\label{Rem_Grigelionis}
Let \ $a \in \RR_+$, \ $b \in \RR$, \ $\sigma \in \RR_{++}$, \ and let \ $m$ \ be a L\'evy
 measure on \ $\RR_{++}$ \ satisfying \eqref{help_Levy}.
Let \ $(Y_t)_{t\in\RR_+}$ \ be the unique strong solution of the SDE
 \eqref{jump_CIR} satisfying \ $\PP(Y_0 \in \RR_+) = 1$ \ and \ $\EE(Y_0) < \infty$.
\ By \eqref{jump_CIR}, the process \ $(Y_t)_{t\in\RR_+}$ \ is a
 semimartingale, see, e.g., Jacod and Shiryaev \cite[I.4.34]{JSh}.
By \eqref{LK}, we have
 \[
   \EE(\ee^{\ii \theta J_t})
   = \exp\left\{\ii t \theta \int_0^1 z \, m(\dd z)
                + t \int_0^\infty
                     \bigl(\ee^{\ii \theta z} - 1 - \ii \theta z h(z)\bigr)
                     \, m(\dd z) \right\}
 \]
 for \ $\theta \in \RR$ \ and \ $t \in \RR_+$, \ where
 \ $h(z) := z \bbone_{[-1,1]}(z)$, \ $z \in \RR $.
\ Using again \eqref{jump_CIR} and the L\'evy--It\^{o}'s representation \eqref{J_Levy_Ito}
 of \ $J$, \ we can write the process \ $(Y_t)_{t\in\RR_+}$ \ in the form
 \begin{equation}\label{Grigelionis}
  \begin{aligned}
   Y_t &= Y_0 + \int_0^t (a - b Y_u) \, \dd u + t \int_0^1 z \, m(\dd z)
          + \sigma \int_0^t \sqrt{Y_u} \, \dd W_u \\
       &\quad
          + \int_0^t \int_\RR
             h(z) \, \tmu^J(\dd u, \dd z)
          + \int_0^t \int_\RR (z - h(z)) \, \mu^J(\dd u, \dd z) , \qquad t \in \RR_+ ,
  \end{aligned}
 \end{equation}
 where \ $\tmu^J(\dd s, \dd z) := \mu^J(\dd s, \dd z) - \dd s\,m(\dd z)$.
\ In fact, \eqref{Grigelionis} is a so-called Grigelionis form for the semimartingale
 \ $(Y_t)_{t\in\RR_+}$, \ see, e.g., Jacod and Shiryaev \cite[III.2.23]{JSh} or
 Jacod and Protter \cite[Theorem 2.1.2]{JacPro}.
\proofend
\end{Rem}

Next we give a statistic for \ $\sigma^2$ \ using continuous time observations
 \ $(Y_t)_{t\in[0,T]}$ \ with some \ $T > 0$.
\ Due to this result we do not consider the estimation of the parameter \ $\sigma$, \ it is
 supposed to be known.

\begin{Rem}\label{Thm_MLE_cons_sigma}
Let \ $a \in \RR_+$, \ $b \in \RR$, \ $\sigma \in \RR_{++}$, \ and let \ $m$ \ be a L\'evy
 measure on \ $\RR_{++}$ \ satisfying \eqref{help_Levy}.
Let \ $(Y_t)_{t\in\RR_+}$ \ be the unique strong solution of the SDE \eqref{jump_CIR}
 satisfying \ $\PP(Y_0 \in \RR_+) = 1$ \ and \ $\EE(Y_0) < \infty$.
\ The Grigelionis representation given in \eqref{Grigelionis} implies that the
 continuous martingale part \ $Y^\cont$ \ of \ $Y$ \ is
 \ $Y^\cont_t = \sigma \int_0^t \sqrt{Y_u} \, \dd W_u$, \ $t \in \RR_+$, \ see Jacod and
 Shiryaev \cite[III.2.28 Remarks, part 1)]{JSh}.
Consequently, the (predictable) quadratic variation process of \ $Y^\cont$ \ is
 \ $\langle Y^\cont \rangle_t = \sigma^2 \int_0^t Y_u \, \dd u$, \ $t \in \RR_+$.
\ Suppose that we have \ $\PP(Y_0 \in \RR_{++}) = 1$ \ or \ $a \in \RR_{++}$.
\ Then for all \ $T \in\RR_{++}$, \ we have
 \[
   \sigma^2 = \frac{\langle Y^\cont \rangle_T}{\int_0^T Y_u \, \dd u} =: \hsigma^2_T ,
 \]
 since, due to Proposition \ref{Pro_jump_CIR},
 \ $\PP\bigl(\int_0^T Y_u \, \dd u \in \RR_{++}\bigr) = 1$.
\ We note that \ $\hsigma^2_T$ \ is a statistic, i.e., there exists a measurable
 function \ $\Xi : D([0,T], \RR) \to \RR$ \ such that
 \ $\hsigma^2_T = \Xi((Y_u)_{u\in[0,T]})$, \ where \ $D([0,T], \RR)$ \ denotes the
 space of real-valued c\`adl\`ag functions defined on \ $[0,T]$, \ since
 \begin{equation}\label{sigma}
  \frac{1}{\frac{1}{n} \sum_{i=1}^{\lfloor nT\rfloor} Y_{\frac{i-1}{n}}}
  \Biggl(\sum_{i=1}^{\lfloor nT\rfloor} \bigl(Y_{\frac{i}{n}} - Y_{\frac{i-1}{n}}\bigr)^2
         - \sum_{u\in[0,T]} (\Delta Y_u)^2\Biggr)
  \stoch \hsigma^2_T \qquad \text{as \ $n \to \infty$,}
 \end{equation}
 where the convergence in \eqref{sigma} holds almost surely along a suitable subsequence,
 the members of the sequence in \eqref{sigma} are measurable functions of
 \ $(Y_u)_{u\in[0,T]}$, \ and one can use Theorems 4.2.2 and 4.2.8 in Dudley \cite{Dud}.
Next we prove \eqref{sigma}.
By Theorem I.4.47 a) in Jacod and Shiryaev \cite{JSh},
 \[
   \sum_{i=1}^{\lfloor nT\rfloor} \bigl(Y_{\frac{i}{n}} - Y_{\frac{i-1}{n}}\bigr)^2
   \stoch [Y]_T \qquad \text{as \ $n \to \infty$,} \qquad
   T \in  \RR_+,
 \]
 where \ $([Y]_t)_{t\in \RR_+}$ \ denotes the quadratic variation process of the
 semimartingale \ $Y$.
\ By Theorem I.4.52 in Jacod and Shiryaev \cite{JSh},
 \[
   [Y]_T = \langle Y^\cont \rangle_T + \sum_{u\in[0,T]} (\Delta Y_u)^2 ,
   \qquad T \in \RR_+ .
 \]
Consequently, for all \ $T \in \RR_+$, \ we have
 \[
   \sum_{i=1}^{\lfloor nT\rfloor} \bigl(Y_{\frac{i}{n}} - Y_{\frac{i-1}{n}}\bigr)^2
   - \sum_{u\in[0,T]} (\Delta Y_u)^2
   \stoch \langle Y^\cont \rangle_T \qquad \text{as \ $n \to \infty$.}
 \]
Moreover, for all \ $T \in \RR_+$, \ we have
 \[
   \frac{1}{n} \sum_{i=1}^{\lfloor nT\rfloor} Y_{\frac{i-1}{n}}
   \stoch \int_0^T Y_u \, \dd u  \qquad \text{as \ $n \to \infty$,}
 \]
 see Proposition I.4.44 in Jacod and Shiryaev \cite{JSh}.
Hence \eqref{sigma} follows by the fact that convergence in probability
 is closed under multiplication.
\proofend
\end{Rem}

\section{Joint Laplace transform of \ $Y_t$ \ and \ $\int_0^t Y_s \, \dd s$}
\label{section_Laplace}

We study the joint Laplace transform of \ $Y_t$ \ and \ $\int_0^t Y_s \, \dd s$, \ since it
 plays a crucial role in deriving the asymptotic behavior of the MLE of \ $b$ \ given in
 \eqref{MLEb}.
Our formula for the joint Laplace transform in question given in Theorem \ref{Thm_Laplace_joint}
 is in accordance with the corresponding one obtained in Keller-Ressel \cite[Theorem 4.10]{KR2}
 in case of a regular affine process and with the one in Jiao et al. \cite[Proposition 4.3]{JiaMaSco}
 in case of a general CBI process.
Here, our contribution is to give a new proof for this joint Laplace transform and to make
 the solutions of the Riccati-type differential equations appearing in the formulas of
 Keller-Ressel \cite[Theorem 4.10]{KR2} and Jiao et al. \cite[Proposition 4.3]{JiaMaSco}
 explicit in case of \ $(Y_t)_{t\in\RR_+}$, \ which turns out to be crucial for our forthcoming statistical study.
For all \ $b \in \RR$ \ and \ $v \in \RR_-$, \ let us introduce the notation
 \ $\gamma_v := \sqrt{b^2 - 2 \sigma^2 v}$.

\begin{Thm}\label{Thm_Laplace_joint}
Let \ $a \in \RR_+$, \ $b \in \RR$, \ $\sigma \in \RR_{++}$, \ and let \ $m$ \ be a L\'evy
 measure on \ $\RR_{++}$ \ satisfying \eqref{help_Levy}.
Let \ $(Y_t)_{t\in\RR_+}$ \ be the unique strong solution of the SDE \eqref{jump_CIR}
 satisfying \ $\PP(Y_0 = y_0) = 1$ \ with some \ $y_0 \in \RR_+$.
\ Then for all \ $u, v \in \RR_-$,
 \[
   \EE\left[\exp\left\{uY_t + v \int_0^t Y_s \, \dd s\right\}\right]
   = \exp\bigg\{\psi_{u,v}(t) y_0
                + \int_0^t
                   \biggl(a \psi_{u,v}(s)
                          + \int_0^\infty
                             \bigl(\ee^{z\psi_{u,v}(s)}
                                   - 1\bigr) m(\dd z)\biggr) \dd s \bigg\}
 \]
 for \ $t \in \RR_+$, \ where the function \ $\psi_{u,v} : \RR_+ \to \RR_-$ \ takes the
 form
 \begin{align}\label{psi_form}
  \psi_{u,v}(t)
  = \begin{cases}
     \frac{u \gamma_v \cosh\left(\frac{\gamma_v t}{2}\right)
           + (-u b + 2 v) \sinh\left(\frac{\gamma_v t}{2}\right)}
          {\gamma_v \cosh\left(\frac{\gamma_v t}{2}\right)
           + (- \sigma^2 u + b) \sinh\left(\frac{\gamma_v t}{2}\right)}
      & \text{if \ $v \in\RR_{--}$ \ or \ $b \ne 0$ \ (i.e., if \ $\gamma_v \in\RR_{++}$),} \\[1mm]
     \frac{u}{1-\frac{\sigma^2u}{2} t}
      & \text{if \ $v = 0$ \ and \ $b = 0$ \ (i.e., if \ $\gamma_v = 0$),}
    \end{cases}
  \qquad t \in \RR_+ .
 \end{align}
\end{Thm}

\begin{Rem}\label{REM_osszhang1}
(i) If \ $v\in\RR_{--}$, \ then \ $\gamma_v > |b|$, \ hence
 \begin{align*}
  \gamma_v\cosh\left(\frac{\gamma_v t}{2}\right)
  + (- \sigma^2 u + b) \sinh\left(\frac{\gamma_v t}{2}\right)
  \geq (\gamma_v + b) \sinh\left(\frac{\gamma_v t}{2}\right)
  \in\RR_{++}, \qquad t\in\RR_+.
 \end{align*}
If \ $v = 0$ \ and \ $b \ne 0$, \ then \ $\gamma_v = |b| \in\RR_{++}$, \ hence
 \begin{align*}
  \gamma_v\cosh\left(\frac{\gamma_v t}{2}\right)
  + (- \sigma^2 u + b) \sinh\left(\frac{\gamma_v t}{2}\right)
  \geq |b| \left(\cosh\left(\frac{\gamma_v t}{2}\right)
                 + \frac{b}{|b|} \sinh\left(\frac{\gamma_v t}{2}\right)\right)
  \in\RR_{++},\qquad t\in\RR_+.
 \end{align*}
Consequently, if \ $\gamma_v \in\RR_{++}$, \ then
 \ $\gamma_v\cosh\left(\frac{\gamma_v t}{2}\right)
    + (- \sigma^2 u + b) \sinh\left(\frac{\gamma_v t}{2}\right) \in\RR_{++}$, $t\in\RR_+$,
 \ hence the function \ $\psi_{u,v}$ \ in \eqref{psi_form} is well-defined.

\noindent (ii)
In Theorem \ref{Thm_Laplace_joint}, we have
 \begin{align}\label{intpsi}
  \int_0^t \psi_{u,v}(s) \, \dd s
  =\begin{cases}
    \frac{b}{\sigma^2} \, t
    -\frac{2}{\sigma^2}
     \log\bigl(\cosh\left(\frac{\gamma_v t}{2}\right)
               + \frac{-\sigma^2 u + b}{\gamma_v}
                 \sinh\left(\frac{\gamma_v t}{2}\right)\bigr) ,
     & \text{if \ $v \in\RR_{--}$ \ or \ $b \ne 0$,} \\
    - \frac{2}{\sigma^2} \log\bigl(1 - \frac{\sigma^2u}{2} t\bigr) ,
     & \text{if \ $v = 0$ \ and \ $b = 0$,}
   \end{cases}
 \end{align}
 for all \ $t \in \RR_+$, \ see, e.g., Lamberton and Lapeyre
 \cite[Chapter 6, Proposition 2.5]{LamLap}.
\proofend
\end{Rem}

\noindent{\bf Proof of Theorem \ref{Thm_Laplace_joint}.}
Introducing the process \ $Z_t := \int_0^t Y_s \, \dd s$, \ $t \in \RR_+$, \ first, we show
 that \ $(Y_t, Z_t)_{t\in\RR_+}$ \ is a 2-dimensional CBI process.
Using the SDE \eqref{jump_CIR} and \eqref{J_Levy_Ito}, this process satisfies a SDE of the
 form given in Barczy et al.\ \cite[Section 5]{BarLiPap2}, namely,
 \begin{align*}
   \begin{bmatrix} Y_t \\ Z_t \end{bmatrix}
   &= \begin{bmatrix} Y_0 \\ 0 \end{bmatrix}
     + \int_0^t
        \Biggl(\begin{bmatrix} a \\ 0 \end{bmatrix}
               + \begin{bmatrix} -b & 0 \\ 1 & 0 \end{bmatrix}
                 \begin{bmatrix} Y_s \\ Z_s \end{bmatrix}\Biggr)
        \dd s
     + \int_0^t \sqrt{\sigma^2 Y_s} \begin{bmatrix}
                                      1 \\
                                      0 \\
                                    \end{bmatrix}
                                    \begin{bmatrix}
                                      1 \\
                                      0 \\
                                    \end{bmatrix}^\top
         \begin{bmatrix} \dd W_s \\ \dd \tW_s \end{bmatrix} \\
   &\quad
     + \int_0^t \sqrt{0\cdot Z_s} \begin{bmatrix}
                                      0 \\
                                      1 \\
                                    \end{bmatrix}
                                    \begin{bmatrix}
                                      0 \\
                                      1 \\
                                    \end{bmatrix}^\top
         \begin{bmatrix} \dd W_s \\ \dd \tW_s \end{bmatrix}
     + \int_0^t \int_{\RR_+^2\setminus\{\bzero\}} \br \, M(\dd s, \dd\br) ,
     \qquad t\in\RR_+,
 \end{align*}
 where \ $(W_t)_{t\in\RR_+}$ \ and \ $(\tW_t)_{t\in\RR_+}$ \ are independent standard Wiener
 processes, and \ $M$ \ is a Poisson random measure on
 \ $\RR_+ \times (\RR_+^2\setminus\{\bzero\})$ \ with intensity measure
 \ $\dd s \, \nu(\dd\br)$, \ where the measure \ $\nu$ \ on \ $\RR_+^2\setminus\{\bzero\}$
 \ is given by \ $\nu(B) := \int_0^\infty \bbone_B(z, 0) \, m(\dd z)$,
 \ $B \in \cB(\RR_+^2\setminus\{\bzero\})$, \ hence
 \ $\int_{\RR_+^2\setminus\{\bzero\}} (1 \land \|\br\|) \, \nu(\dd\br)
    \leq \int_0^\infty z \, m(\dd z) < \infty$.
\ Put
 \[
   \bc := \begin{bmatrix} c_1 \\ c_2 \end{bmatrix}
   := \begin{bmatrix} \frac{1}{2} \, \sigma^2 \\ 0 \end{bmatrix} \in \RR_+^2 , \qquad
   \bbeta := \begin{bmatrix} a \\ 0 \end{bmatrix} \in \RR_+^2 , \qquad
   \bB := \begin{bmatrix} -b & 0 \\ 1 & 0 \end{bmatrix} ,\qquad
   \bmu := (\mu_1, \mu_2) := (0, 0) .
 \]
Then \ $\bB$ \ is an essentially non-negative matrix (i.e., its off-diagonal entries are
 non-negative), and, due to \eqref{help_Levy},
 \[
   \int_{\RR_+^2\setminus\{\bzero\}} \|\br\| \bbone_{\{ \|\br\| \geq 1\}} \, \nu(\dd\br)
   \leq \int_{\RR_+^2\setminus\{\bzero\}} \|\br\| \, \nu(\dd\br)
   = \int_0^\infty z \, m(\dd z) < \infty
 \]
 yielding that condition (2.7) of Barczy et al.\ \cite{BarLiPap2} is satisfied.
Thus, by Theorem 4.6 in Barczy et al.\ \cite{BarLiPap2}, \ $(Y_t, Z_t)_{t\in\RR_+}$ \ is a
 CBI process with parameters \ $(2, \bc, \bbeta, \bB, \nu, \bmu)$.
\ We note that the fact that \ $(Y_t, Z_t)_{t\in\RR_+}$ \ is a 2-dimensional CBI process
 follows from Filipovi\'{c} et al. \cite[paragraph before Theorem 4.3]{FilMaySch} as well,
 where this property is stated for general affine processes without a proof.
The branching mechanism of \ $(Y_t, Z_t)_{t\in\RR_+}$ \ is
 \ $\bR(u, v) = (R_1(u, v), R_2(u, v))$, \ $u, v \in \RR_-$, \ with
 \[
   R_1(u, v)
   = c_1 u^2
     + \left\langle \bB \begin{bmatrix} 1 \\ 0 \end{bmatrix},
                    \begin{bmatrix} u \\ v \end{bmatrix}\right\rangle
   = \frac{\sigma^2}{2} u^2 - b u +  v , \qquad
   R_2(u, v)
   = c_2 v^2 + \left\langle \bB \begin{bmatrix} 0 \\ 1 \end{bmatrix},
                    \begin{bmatrix} u \\ v \end{bmatrix}\right\rangle
   = 0 ,
 \]
 and the immigration mechanism of \ $(Y_t, Z_t)_{t\in\RR_+}$ \ is
 \[
   F(u, v)
   = \left\langle \bbeta, \begin{bmatrix} u \\ v \end{bmatrix}\right\rangle
     - \int_{\RR_+^2\setminus\{\bzero\}}
        \Biggl(\exp\Biggl\{\left\langle\begin{bmatrix} u \\ v \end{bmatrix},
                                       \br\right\rangle\Biggr\}
         - 1\Biggr)
        \, \nu(\dd\br)
   = a u + \int_0^\infty (\ee^{uz} - 1) \, m(\dd z) ,
 \]
 see, e.g., Theorem 2.4 in Barczy et al.\ \cite{BarLiPap2}.
Note that \ $\bR(u,v) = (R(u) +v,0)$, $u,v\in\RR_{-}$, \ and
 \ $F(u,v) = F(u)$, $u,v\in\RR_{-}$, \ where \ $R(u)$, $u\in\RR_{-}$, \ and \ $F(u)$, $u\in\RR_{-}$, \ are given
 in Proposition \ref{Pro_jump_CIR}, which are in accordance with Theorem 4.10 in Keller-Ressel \cite{KR2}.
Consequently, by Theorem 2.7 of Duffie et al.\ \cite{DufFilSch}
 (see also Barczy et al.\ \cite[Theorem 2.4]{BarLiPap2}), we have
 \begin{equation}\label{Laplace_joint}
  \EE\left[\exp\left\{uY_t + v \int_0^t Y_s \, \dd s\right\}\right]
  = \exp\bigg\{\psi_{u,v}(t) y_0
               + \int_0^t F(\psi_{u,v}(s), \varphi_{u,v}(s)) \, \dd s \bigg\}
 \end{equation}
 for \ $t \in \RR_+$, \ where the function
 \ $(\psi_{u,v}, \varphi_{u,v}) : \RR_+ \to \RR_-^2$ \ is the unique locally bounded
 solution to the system of differential equations
 \[
   \begin{cases}
    \psi'_{u,v}(t)
    = R_1(\psi_{u,v}(t), \varphi_{u,v}(t))
    = \frac{\sigma^2}{2} \psi_{u,v}(t)^2 - b \psi_{u,v}(t) + \varphi_{u,v}(t) ,
    \qquad t \in \RR_+ , \\
   \varphi'_{u,v}(t)
   = R_2(\psi_{u,v}(t), \varphi_{u,v}(t))
   = 0 , \qquad t \in \RR_+ ,
  \end{cases}
 \]
 with initial values \ $\psi_{u,v}(0) = u$, \ $\varphi_{u,v}(0) = v$.
\ Clearly, \ $\varphi_{u,v}(t) = v$, \ $t \in \RR_+$, \ hence we obtain
 \[
   \psi'_{u,v}(t) = \frac{\sigma^2}{2} \psi_{u,v}(t)^2 - b \psi_{u,v}(t) + v ,
   \qquad t \in \RR_+ , \qquad
   \psi_{u,v}(0) = u .
 \]
The solution of this differential equation is \eqref{psi_form}.
Indeed, in case of \ $\gamma_v > 0$, \ we can refer to, e.g., Lamberton and Lapeyre
 \cite[Chapter 6, Proposition 2.5]{LamLap}, and in case of \ $\gamma_v = 0$, \ this is a simple
 separable ODE taking the form \ $\psi_{u,0}'(t) = \frac{\sigma^2}{2} \psi_{u,0}(t)^2$,
 \ $t \in \RR_+$, \ with initial condition \ $\psi_{u,0}(0) = u$.
\ Hence, by \eqref{Laplace_joint}, we obtain the statement.
\proofend

\begin{Ex}\label{Ex_Laplace_joint_crit}
Now we formulate a special case of Theorem \ref{Thm_Laplace_joint} in the critical case
 \ ($b=0$) \ supposing that the L\'evy measure \ $m$ \ takes the form given in \eqref{Levy_cond},
 i.e., in the case of a critical BAJD process.
For all \ $u, v \in\RR_-$, \ let us introduce the notations
 \begin{align*}
  \alpha^{(1)}_{u,v} := u \gamma_v  + 2v , \quad
  \alpha^{(2)}_{u,v} := u \gamma_v  - 2v , \quad
  \beta^{(1)}_{u,v} := \lambda (-\sigma^2 u + \gamma_v) - \alpha^{(1)}_{u,v} , \quad
  \beta^{(2)}_{u,v} := \lambda (\sigma^2 u + \gamma_v) - \alpha^{(2)}_{u,v} ,
 \end{align*}
 where \ $\gamma_v=\sqrt{-2\sigma^2 v}$ \ (since now \ $b = 0$).
\ Let \ $(Y_t)_{t\in\RR_+}$ \ be the unique strong solution of the SDE \eqref{jump_CIR}
 satisfying \ $\PP(Y_0 = y_0) = 1$ \ with some \ $y_0 \in \RR_+$, \ with \ $b=0$ \ and \ $m$ \ being
 a L\'evy measure on \ $\RR_{++}$ \ satisfying \eqref{Levy_cond}.
Then we check that for all \ $u, v \in \RR_-$,
 \begin{align}\label{help_BAJD_Laplace1}
   \EE\left[\exp\left\{uY_t + v \int_0^t Y_s \, \dd s\right\}\right]
   = \exp\big\{\psi_{u,v}(t) y_0 + \phi_{u,v}(t) \big\} ,
   \qquad t \in \RR_+ ,
 \end{align}
 where the function \ $\psi_{u,v} : \RR_+ \to \RR_-$ \ is given by \eqref{psi_form} with \ $\gamma_v=\sqrt{-2\sigma^2 v}$ \
  (since now \ $b = 0$), \ and if \ $v \in\RR_{--}$ \ (i.e., if \ $\gamma_v \in\RR_{++}$) \ and
 \ $v \notin \{-\frac{\sigma^2u^2}{2},-\frac{\sigma^2\lambda^2}{2}\}$
 \ (i.e., if \ $\beta^{(2)}_{u,v} \ne 0$), \ then
 \begin{align}\label{phi_1}
  \begin{split}
   \phi_{u,v}(t)
   &= - \frac{2a}{\sigma^2}\log\left(\cosh\left(\frac{\gamma_v t}{2}\right)
      - \frac{\sigma^2 u}{\gamma_v} \sinh\left(\frac{\gamma_v t}{2}\right)\right) \\
   &\quad
      + c \left(\frac{\alpha^{(2)}_{u,v} }{\beta^{(2)}_{u,v}} t
                + \frac{1}{\gamma_v}\left( \frac{\alpha^{(1)}_{u,v}}{\beta^{(1)}_{u,v}}
                - \frac{\alpha^{(2)}_{u,v}}{\beta^{(2)}_{u,v}}\right)
                  \log\left(\frac{\beta^{(1)}_{u,v}\ee^{\gamma_v t}
                            + \beta^{(2)}_{u,v}}{\beta^{(1)}_{u,v}+\beta^{(2)}_{u,v}}\right)
          \right) , \qquad t \in \RR_+ ,
  \end{split}
 \end{align}
 and if \ $v \in \{-\frac{\sigma^2u^2}{2}, -\frac{\sigma^2\lambda^2}{2}\}$ \ and \ $v \in\RR_{--}$
 \ (i.e., if \ $v \in\RR_{--}$ \ and \ $\beta^{(2)}_{u,v} = 0$), \ then
 \begin{align}\label{phi_3}
  \begin{split}
   \phi_{u,v}(t)
   & = - \frac{2a}{\sigma^2}\log\left( \cosh\left(\frac{\gamma_v t}{2}\right)
       - \frac{\sigma^2 u}{\gamma_v} \sinh\left(\frac{\gamma_v t}{2}\right)\right) \\
   &\quad
       + \frac{c}{\beta^{(1)}_{u,v}}
         \left(\alpha^{(1)}_{u,v} t
               + \frac{\alpha^{(2)}_{u,v}}{\gamma_v}(1- \ee^{-\gamma_v t})\right) ,
   \qquad t \in \RR_+ ,
  \end{split}
 \end{align}
 and if \ $v = 0$ \ (i.e., if \ $\gamma_v = 0$), \ then
 \begin{align}\label{phi_2}
  \phi_{u,v}(t)
  = -\frac{2a}{\sigma^2}\log\left(1-\frac{\sigma^2 u}{2}t\right)
    - \frac{2c}{\sigma^2 \lambda}
      \log\left(1 - \frac{\sigma^2 \lambda u}{2(\lambda-u)}t\right) ,
  \qquad t \in \RR_+ .
 \end{align}
Especially, for all \ $u \in \RR_-$,
 \[
   \EE(\ee^{uY_t})
   = \exp\left\{\frac{uy_0}{1-\frac{\sigma^2 u}{2}t}\right\}
     \left(1 - \frac{\sigma^2 u}{2}t\right)^{-\frac{2a}{\sigma^2}}
     \left(1 - \frac{\sigma^2 \lambda u}{2(\lambda-u)}t\right)^{-\frac{2c}{\sigma^2 \lambda}},
   \qquad t \in \RR_+ .
 \]

First, we check that the function \ $\phi_{u,v}$ \ in \eqref{phi_1} is well-defined.
If \ $\gamma_v \in \RR_{++}$, \ then \ $\beta^{(1)}_{u,v} \in \RR_{++}$ \ (due to
 \ $\alpha^{(1)}_{u,v}\in\RR_-$).
\ If \ $\gamma_v = 0$, \ then \ $\alpha^{(1)}_{u,v} = \alpha^{(2)}_{u,v} = 0$,
 \ $\beta^{(1)}_{u,v} = -\sigma^2 \lambda u$ \ and
 \ $\beta^{(2)}_{u,v} = \sigma^2 \lambda u$.
\ Further,
 \ $\beta^{(1)}_{u,v} \ee^{\gamma_v t} + \beta^{(2)}_{u,v}
    \geq \beta^{(1)}_{u,v} + \beta^{(2)}_{u,v} \in \RR_{+}$,
 \ and \ $\beta^{(1)}_{u,v} + \beta^{(2)}_{u,v} = 0$ \ holds if and only if \ $\gamma_v = 0$
 \ (i.e., \ $v = 0$).
\ Indeed, since \ $\beta^{(1)}_{u,v} \in \RR_+$, \ we have
 \[
   \beta^{(1)}_{u,v} \ee^{\gamma_v t} + \beta^{(2)}_{u,v}
   \geq \beta^{(1)}_{u,v} + \beta^{(2)}_{u,v} = 2\gamma_v(\lambda-u) \in \RR_+ ,
 \]
 and \ $\beta^{(1)}_{u,v} + \beta^{(2)}_{u,v} = 0$ \ holds if and only if \ $\gamma_v = 0$
 \ (since \ $\lambda - u \in \RR_{++}$ \ due to \ $\lambda \in \RR_{++}$ \ and
 \ $u \in \RR_{-}$).
\ Consequently, if \ $\gamma_v \in \RR_{++}$, \ then \ $\beta^{(1)}_{u,v} \in \RR_{++}$ \ and
 \ $\beta^{(1)}_{u,v} + \beta^{(2)}_{u,v} \in \RR_{++}$, \ yielding that the function
 \ $\phi_{u,v}$ \ in \eqref{phi_1} is well-defined.

Next, we check that, provided that \ $v \in \RR_{--}$ \ (i.e., \ $\gamma_v \in \RR_{++}$), \ we have
 \ $\beta^{(2)}_{u,v} = 0$ \ if and only if
 \ $v \in \{-\frac{\sigma^2u^2}{2}, -\frac{\sigma^2\lambda^2}{2}\}$.
\ Indeed, \ $\beta^{(2)}_{u,v} = 0$ \ holds if and only if
 \ $\lambda(\sigma^2 u + \gamma_v) = \gamma_v u - 2v$, \ i.e., by some algebraic
 transformations, if and only if
 \[
   (\sqrt{-2v})^2 + \sigma(u - \lambda) \sqrt{-2v} - \lambda \sigma^2 u = 0 .
 \]
Solving this quadratic equation with respect to \ $\sqrt{-2v}$, \ we get
 \ $\sqrt{-2v} \in \{-\sigma u, \sigma \lambda\}$, \ yielding that
 \ $v \in \{-\frac{\sigma^2u^2}{2}, -\frac{\sigma^2\lambda^2}{2}\}$, \ as desired.

Finally, we check \eqref{help_BAJD_Laplace1}.
Using Theorem \ref{Thm_Laplace_joint} and part (ii) of Remark \ref{REM_osszhang1}, it
 remains to calculate
 \ $\int_0^t \left(\int_0^\infty (\ee^{z\psi_{u,v}(s)} - 1) \, m(\dd z)\right)\dd s$
 \ for \ $u, v \in \RR_{-}$.
\ For all \ $u, v \in \RR_{-}$, \ we have
 \begin{align}\label{help_crit_Lap_spec}
  \begin{aligned}
   &\int_0^\infty (\ee^{z\psi_{u,v}(s)} - 1) \, m(\dd z)
    = \int_0^\infty (\ee^{z\psi_{u,v}(s)} - 1) c \lambda \ee^{-\lambda z} \, \dd z
    = c \lambda \int_0^\infty \ee^{(\psi_{u,v}(s)-\lambda)z} \, \dd z - c \\
   &= \frac{c\lambda}{\psi_{u,v}(s) - \lambda}
      \left(\lim_{z\to\infty} \ee^{(\psi_{u,v}(s)-\lambda)z}  -1\right) - c
    = -\frac{c\psi_{u,v}(s)}{\psi_{u,v}(s) - \lambda} ,
    \qquad s \in \RR_+ ,
  \end{aligned}
 \end{align}
 where we used that \ $\psi_{u,v}(s) - \lambda \in \RR_{--}$, \ $s \in \RR_+$, \ following
 from \ $\lambda \in \RR_{++}$ \ and \ $\psi_{u,v}(s) \in \RR_{-}$, \ $s \in \RR_+$.
\ By some algebraic transformations,
 \begin{align*}
  \frac{c\psi_{u,v}(t)}{\lambda-\psi_{u,v}(t)}
  &= c \frac{u\gamma_v(1+\ee^{\gamma_v t})+2v(\ee^{\gamma_v t}-1)}
            {-\lambda\sigma^2 u(\ee^{\gamma_v t}-1)+\lambda\gamma_v
             +\lambda\gamma_v\ee^{\gamma_v t}
             -u\gamma_v(1+\ee^{\gamma_v t})-2v(\ee^{\gamma_v t}-1)} \\
  &= c \frac{\alpha^{(1)}_{u,v}\ee^{\gamma_v t}+\alpha^{(2)}_{u,v}}
            {\beta^{(1)}_{u,v}\ee^{\gamma_v t}+\beta^{(2)}_{u,v}} ,
  \qquad t \in \RR_+ .
 \end{align*}
As we have seen, if \ $v \in \RR_{--}$ \ (i.e., if
 \ $\gamma_v \in \RR_{++}$), \ then \ $\beta^{(1)}_{u,v} \in \RR_{++}$,
 \ and \ $\beta^{(2)}_{u,v} = 0$ \ if and only if
 \ $v \in \{-\frac{\sigma^2u^2}{2}, -\frac{\sigma^2\lambda^2}{2}\}$.
\ Hence
 \[
   \frac{c\psi_{u,v}(t)}{\lambda-\psi_{u,v}(t)}
   = \begin{cases}
      c\left(\frac{\alpha^{(2)}_{u,v}}{\beta^{(2)}_{u,v}}
             + \left(\alpha^{(1)}_{u,v}
                     - \frac{\beta^{(1)}_{u,v}\alpha^{(2)}_{u,v}}{\beta^{(2)}_{u,v}}\right)
               \frac{\ee^{\gamma_v t}}
                    {\beta^{(1)}_{u,v}\ee^{\gamma_v t}+\beta^{(2)}_{u,v}}\right)
       & \text{if \ $v \in \RR_{--}$,
               \ $v \notin \{-\frac{\sigma^2u^2}{2}, -\frac{\sigma^2 \lambda^2}{2}\}$,} \\
      \frac{c}{\beta^{(1)}_{u,v}} (\alpha^{(1)}_{u,v} + \alpha^{(2)}_{u,v} \ee^{-\gamma_v t})
       & \text{if \ $v \in \{-\frac{\sigma^2u^2}{2}, -\frac{\sigma^2 \lambda^2}{2}\}$,
               \ $v \in \RR_{--}$,} \\
      \frac{2cu}{2(\lambda-u) - \lambda u\sigma^2t}
       & \text{if \ $v = 0$,}
      \end{cases}
 \]
 for \ $t \in \RR_+$.
\ By integration, we get \eqref{phi_1}, \eqref{phi_3} and \eqref{phi_2}.
\proofend
\end{Ex}

Next, we formulate two corollaries of Theorem \ref{Thm_Laplace_joint} giving the Laplace
 transforms of \ $Y_t$ \ and \ $\int_0^t Y_s \, \dd s$, \ $t \in \RR_+$, \ separately.

\begin{Cor}\label{Thm_Laplace_Y}
Let \ $a \in \RR_+$, \ $b \in \RR$, \ $\sigma \in \RR_{++}$, \ and let \ $m$ \ be a L\'evy
 measure on \ $\RR_{++}$ \ satisfying \eqref{help_Levy}.
Let \ $(Y_t)_{t\in\RR_+}$ \ be the unique strong solution of the SDE \eqref{jump_CIR}
 satisfying \ $\PP(Y_0 = y_0) = 1$ \ with some \ $y_0 \in \RR_+$.
\ Then for all \ $u \in \RR_-$,
 \[
   \EE\bigl(\ee^{uY_t}\bigr)
   = \exp\biggl\{\psi_{u,0}(t) y_0
                 + \int_0^t
                    \biggl(a \psi_{u,0}(s)
                           + \int_0^\infty
                              \bigl(\ee^{z\psi_{u,0}(s)}
                                    - 1\bigr) m(\dd z)\biggr) \dd s \bigg\} ,
  \qquad t \in \RR_+ ,
 \]
 where the function \ $\psi_{u,0} : \RR_+ \to \RR_-$ \ takes the form
 \begin{align}\label{psi_form_Y}
  \psi_{u,0}(t)
  = \begin{cases}
     \frac{2ub\ee^{-bt}}{\sigma^2u(\ee^{-bt}-1)+2b} & \text{if \ $b \ne 0$,} \\
     \frac{u}{1-\frac{\sigma^2 u}{2}\,t} & \text{if \ $b = 0$,}
    \end{cases}
  \qquad t \in \RR_+ .
 \end{align}
\end{Cor}

\noindent{\bf Proof.}
We have to substitute \ $v = 0$ \ in \eqref{psi_form}.
The case of \ $b = 0$ \ is trivial.
In case of \ $b \ne 0$ \ we have \ $\gamma_v=\gamma_0 = |b|$, \ hence
 \[
   \psi_{u,0}(t)
   = \frac{u |b| \cosh\bigl(\frac{|b| t}{2}\bigr)
           - u b \sinh\bigl(\frac{|b| t}{2}\bigr)}
          {|b| \cosh\bigl(\frac{|b| t}{2}\bigr)
           + (- \sigma^2 u + b) \sinh\bigl(\frac{|b| t}{2}\bigr)} ,
           \qquad t\in\RR_+,
 \]
 yielding \eqref{psi_form_Y}.
\proofend

\begin{Ex}
Now we formulate a special case of Corollary \ref{Thm_Laplace_Y} in the supercritical case
 \ ($b\in\RR_{--}$) \ supposing that the L\'evy measure \ $m$ \ takes the form given in \eqref{Levy_cond},
 i.e., in the case of a supercritical BAJD process.
\ Let \ $(Y_t)_{t\in\RR_+}$ \ be the unique strong solution of the SDE \eqref{jump_CIR}
 satisfying \ $\PP(Y_0 = y_0) = 1$ \ with some \ $y_0 \in \RR_+$, \ with \ $b\in\RR_{--}$ \ and \ $m$ \ being
 a L\'evy measure on \ $\RR_{++}$ \ satisfying \eqref{Levy_cond}.
Then for all \ $u \in \RR_-$,
 \begin{align}\label{help_BAJD_Laplace2}
   \EE(\ee^{uY_t}) = \exp\Big\{\psi_{u,0}(t) y_0 + \phi_{u,0}(t) \Big\} , \qquad t \in \RR_+ ,
 \end{align}
 where the function \ $\psi_{u,0} : \RR_+ \to \RR_-$ \ is given by \eqref{psi_form_Y}, and
 if \ $u \ne \frac{2b}{\sigma^2}$, \ then
 \begin{align}\label{phi_form1_super}
  \phi_{u,0}(t)
  = - \frac{2a}{\sigma^2}\log\left( 1 + \frac{\sigma^2 u}{2b}(\ee^{-bt} - 1) \right)
    + \frac{2c}{-\sigma^2 \lambda + 2b}
      \log\left(\frac{(-\sigma^2\lambda+2b)u\ee^{-bt}+(\sigma^2 u-2b)\lambda}
                     {2b(u-\lambda)}\right)
 \end{align}
 for \ $t \in \RR_+$, \ and if \ $u = \frac{2b}{\sigma^2}$, \ then
 \begin{align}\label{phi_form2_super}
  \phi_{u,0}(t) = \frac{2b(2ab-c\sigma^2-a\sigma^2\lambda)}{\sigma^2(-\sigma^2\lambda+2b)} t ,
  \qquad t \in \RR_+ .
 \end{align}

To check \eqref{help_BAJD_Laplace2}, using Corollary \ref{Thm_Laplace_Y}, it remains to
 calculate
 \ $\int_0^t
     \bigl(a \psi_{u,0}(s) + \int_0^\infty \bigl(\ee^{z\psi_{u,0}(s)} - 1\bigr) m(\dd z)\bigr)
     \dd s$.
\ Here, by \eqref{psi_form_Y},
 \begin{align*}
  \int_0^t \psi_{u,0}(s)\,\dd s
  = -\frac{2}{\sigma^2} \log\left(1 + \frac{\sigma^2 u}{2b}(\ee^{-bt} - 1)\right) ,
  \qquad t \in \RR_+ .
 \end{align*}
Similarly to \eqref{help_crit_Lap_spec}, for all \ $u \in \RR_{-}$,
 \begin{align}
  \int_0^\infty (\ee^{z\psi_{u,0}(s)} - 1) \, m(\dd z)
  = - \frac{c\psi_{u,0}(s)}{\psi_{u,0}(s) - \lambda} ,
  \qquad s \in \RR_+ ,
 \end{align}
By some algebraic transformations,
 \begin{align*}
  - \frac{c\psi_{u,0}(s)}{\psi_{u,0}(s) - \lambda}
  = \frac{-2bcu\ee^{-bt}}{(-\sigma^2 \lambda + 2b)u\ee^{-bt} + (\sigma^2 u - 2b)\lambda} ,
  \qquad t \in \RR_+ ,
 \end{align*}
 where
 \ $(-\sigma^2 \lambda + 2b) u \ee^{-bt} + \sigma^2 \lambda u - 2b\lambda
    \geq (-\sigma^2 \lambda + 2b)u + \sigma^2 \lambda u - 2b \lambda
    = 2 b (u - \lambda) > 0$,
 \ since \ $-\sigma^2 \lambda + 2 b < 0$ \ due to \ $b < 0$ \ and \ $\lambda > 0$.
\ If \ $u = \frac{2b}{\sigma^2}$ \ (i.e., if \ $\sigma^2 u - 2b = 0$), \ then
 \[
   - \int_0^t \frac{c\psi_{u,0}(s)}{\psi_{u,0}(s) - \lambda} \, \dd s
   = \frac{-2bc}{-\sigma^2 \lambda + 2b} t , \qquad t \in \RR_+ ,
 \]
 and hence
 \begin{align*}
  \phi_{u,0}(t)
  &= -\frac{2}{\sigma^2} \log\left(1 + \frac{\sigma^2 u}{2b}(\ee^{-bt} - 1)\right)
     - \frac{2bc}{-\sigma^2 \lambda + 2b} t
   =  \frac{2ab}{\sigma^2}t - \frac{2bc}{-\sigma^2 \lambda + 2b} t \\
  &= \frac{2b(2ab - c\sigma^2 - a\sigma^2 \lambda)}{\sigma^2(-\sigma^2 \lambda + 2b)} t ,
  \qquad t \in \RR_+ ,
 \end{align*}
 yielding \eqref{phi_form2_super}.
If \ $u \ne \frac{2b}{\sigma^2}$ \ (i.e., if \ $\sigma^2 u - 2 b \ne 0$), \ then
 \[
   - \int_0^t \frac{c\psi_{u,0}(s)}{\psi_{u,0}(s)-\lambda} \, \dd s
   = \frac{2c}{-\sigma^2 \lambda + 2b}
     \log\left(\frac{(-\sigma^2\lambda+2b)u\ee^{-bt}+(\sigma^2 u-2b)\lambda}
                   {2b(u-\lambda)}\right) ,
   \qquad t \in \RR_+ ,
 \]
 yielding \eqref{phi_form1_super}.
\proofend
\end{Ex}

Substituting \ $u = 0$ \ in \eqref{psi_form}, we obtain the following corollary.

\begin{Cor}\label{Thm_Laplace_intY}
Let \ $a \in \RR_+$, \ $b \in \RR$, \ $\sigma \in \RR_{++}$, \ and let \ $m$ \ be a L\'evy
 measure on \ $\RR_{++}$ \ satisfying \eqref{help_Levy}.
Let \ $(Y_t)_{t\in\RR_+}$ \ be the unique strong solution of the SDE \eqref{jump_CIR}
 satisfying \ $\PP(Y_0 = y_0) = 1$ \ with some \ $y_0 \in \RR_+$.
\ Then for all \ $v \in \RR_-$,
 \[
   \EE\left[\exp\left\{v \int_0^t Y_s \, \dd s\right\}\right]
   = \exp\bigg\{\psi_{0,v}(t) y_0
                + \int_0^t
                   \biggl(a \psi_{0,v}(s)
                          + \int_0^\infty
                             \bigl(\ee^{z\psi_{0,v}(s)}
                                   - 1\bigr) m(\dd z)\biggr) \dd s \bigg\}
 \]
 where the function \ $\psi_{0,v} : \RR_+ \to \RR_-$ \ takes the form
 \begin{align}\label{psi_form_intY}
  \psi_{0,v}(t)
  = \begin{cases}
     \frac{2 v \sinh\left(\frac{\gamma_v t}{2}\right)}
          {\gamma_v \cosh\left(\frac{\gamma_v t}{2}\right)
           + b \sinh\left(\frac{\gamma_v t}{2}\right)}
      & \text{if \ $v \in\RR_{--}$ \ or \ $b \ne 0$,} \\[1mm]
     0 & \text{if \ $v = 0$ \ and \ $b = 0$,}
    \end{cases}
  \qquad t \in \RR_+ .
 \end{align}
\end{Cor}

\section{Existence and uniqueness of MLE}
\label{section_EUMLE}

In this section, we will consider the jump-type CIR model \eqref{jump_CIR} with known
 \ $a \in \RR_+$, \ $\sigma \in \RR_{++}$, \ L\'evy measure \ $m$ \ satisfying
 \eqref{help_Levy}, and a known deterministic initial value \ $Y_0 = y_0 \in \RR_+$, \ and
 we will consider \ $b \in \RR$ \ as an unknown parameter.

Let \ $\PP_b$ \ denote the probability measure induced by \ $(Y_t)_{t\in\RR_+}$ \ on
 the measurable space \ $(D(\RR_+, \RR), \cD(\RR_+, \RR))$ \ endowed with the
 natural filtration \ $(\cD_t(\RR_+, \RR))_{t\in\RR_+}$, \ see Appendix
 \ref{App_LR}.
Further, for all \ $T \in \RR_{++}$, \ let
 \ $\PP_{b,T} := \PP_{b\!\!}|_{\cD_T(\RR_+, \RR)}$ \ be the
 restriction of \ $\PP_b$ \ to \ $\cD_T(\RR_+, \RR)$.

The next proposition is about the form of the Radon--Nikodym derivative
 \ $\frac{\dd \PP_{b,T}}{\dd \PP_{\tb,T}}$ \ for \ $b, \tb \in \RR$.
\ We will consider \ $\PP_{\tb,T}$ \ as a fixed reference measure, and we will derive
 MLE for the parameter \ $b$ \ based on the observations \ $(Y_t)_{t\in[0,T]}$.

\begin{Pro}\label{RNb}
Let \ $b, \tb \in \RR$.
\ Then for all \ $T \in \RR_{++}$, \ the probability measures \ $\PP_{b,T}$ \ and
 \ $\PP_{\tb,T}$ \ are absolutely continuous with respect to each other, and, under \ $\PP$,
 \begin{align}\label{RNformulab}
\log\biggl(\frac{\dd \PP_{b,T}}{\dd \PP_{\tb,T}}(\tY)\biggr)
  = -\frac{b-\tb}{\sigma^2} (\tY_T - y_0 - a T - J_T)
    -\frac{b^2-\tb^2}{2\sigma^2} \int_0^T \tY_s \, \dd s,
 \end{align}
 where \ $\tY$ \ is the process corresponding to the parameter \ $\tb$.
\end{Pro}

\noindent{\bf Proof.}
In what follows, we will apply Theorem III.5.34 in Jacod and Shiryaev \cite{JSh}
 (see also Appendix \ref{App_LR}).
We will work on the canonical space \ $(D(\RR_+, \RR), \cD(\RR_+, \RR))$.
\ Let \ $(\eta_t)_{t\in\RR_+}$ \ denote the canonical process
 \ $\eta_t(\omega) := \omega(t)$, \ $\omega \in D(\RR_+, \RR)$,
 \ $t \in \RR_+$.
\ The jump-type CIR process \eqref{jump_CIR} can be written in the form \eqref{Grigelionis}.
By Proposition \ref{Pro_jump_CIR}, the SDE \eqref{jump_CIR} has a pathwise unique strong
 solution (with the given deterministic initial value \ $y_0 \in \RR_+$), \ and hence, by
 Theorem III.2.26 in Jacod and Shiryaev \cite{JSh}, under the probability measure \ $\PP_b$,
 \ the canonical process \ $(\eta_t)_{t\in\RR_+}$ \ is a semimartingale with semimartingale
 characteristics \ $(B^{(b)}, C, \nu)$ \ associated with the truncation function \ $h$,
 \ where
 \begin{align*}
  &B^{(b)}_t
   =  \int_0^t \left(a - b \eta_u + \int_0^1 z \, m(\dd z)\right) \dd u , \qquad t\in\RR_+,\\
  &C_t = \int_0^t (\sigma \sqrt{\eta_u})^2 \, \dd u
       = \sigma^2 \int_0^t \eta_u \, \dd u , \qquad t\in\RR_+,\\
  &\nu(\dd t, \dd y) = K(\eta_t, \dd y) \, \dd t = \dd t \, m(\dd y)
 \end{align*}
 with the Borel transition kernel \ $K$ \ from \ $\RR_+ \times \RR$ \ into \ $\RR$ \ given
 by
 \[
   K(y, R) := \int_{\RR} \!\bbone_{R\setminus\{0\}}(z) \, m(\dd z) = m(R)
   \qquad \text{for \ $y \in \RR_+$ \ and \ $R \in \cB(\RR)$.}
 \]
The aim of the following discussion is to check the set of sufficient conditions presented
 in Appendix \ref{App_LR} (of which the notations will be used) in order to have right to
 apply Theorem III.5.34 in Jacod and Shiryaev \cite{JSh}.
First note that \ $(C_t)_{t\in\RR_+}$ \ and \ $\nu(\dd t, \dd y)$ \ do not depend on the
 unknown parameter \ $b$, \ and hence \ $V^{(\tb,b)}$ \ is identically one and then
 \eqref{GIR1} and \eqref{GIR2} readily hold.
We also have
 \[
   \PP_b\big(\nu( \{t\}\times \RR) = 0 \big) = \PP_b(0\cdot m(\RR) =0) = 1 ,
   \qquad t \in \RR_+ , \quad b \in \RR ,
 \]
 since \ $m(\RR) = m(\RR_+) \in\RR_+$ \ due to \eqref{help_Levy}.
Further, \ $(C_t)_{t\in\RR_+}$ \ can be represented as \ $C_t = \int_0^t c_u \, \dd F_u$,
 \ $t \in \RR_+$, \ where the stochastic processes \ $(c_t)_{t\in\RR_+}$ \ and
 \ $(F_t)_{t\in\RR_+}$ \ are given by \ $c_t := \sigma^2 \eta_t$, \ $t \in \RR_+$, \ and
 \ $F_t = t$, \ $t \in \RR_+$.
\ Consequently, for all \ $b, \tb \in \RR$,
 \begin{align*}
  B^{(b)}_t - B^{(\tb)}_t
  = - (b - \tb) \int_0^t \eta_u \, \dd u
  = \int_0^t c_u \beta^{(\tb,b)}_u \, \dd F_u
 \end{align*}
 $\PP_b$-almost surely for every \ $t \in \RR_+$, \ where the stochastic process
 \ $(\beta^{(\tb,b)}_t)_{t\in\RR_+}$ \ is given by
 \[
   \beta^{(\tb,b)}_t = - \frac{b - \tb}{\sigma^2} , \qquad t \in \RR_+ ,
 \]
 which yields \eqref{GIR3}.

Next we check \eqref{GIR4}, i.e.,
 \begin{gather}\label{COND}
  \PP_b\left(\int_0^t \bigl(\beta^{(\tb,b)}_u\bigr)^2 c_u \, \dd F_u < \infty\right)
  = 1 ,  \qquad t \in \RR_+ .
 \end{gather}
We have
 \[
   \int_0^t \bigl(\beta^{(\tb,b)}_u\bigr)^2 c_u \, \dd F_u
   = \frac{(b-\tb)^2}{\sigma^2} \int_0^t \eta_u \, \dd u , \qquad t \in \RR_+ .
 \]
Since for each \ $\omega \in D(\RR_+, \RR)$, \ the trajectory
 \ $[0, t] \ni u \mapsto \eta_u(\omega)$ \ is c\`adl\`ag, hence bounded (see, e.g.,
 Billingsley \cite[(12.5)]{Bil}), we have \ $\int_0^t \eta_u(\omega) \, \dd u < \infty$,
 \ hence we obtain \eqref{COND}.

Next, we check that, under the probability measure \ $\PP_{b}$, \ local uniqueness holds for
 the martingale problem on the canonical space corresponding to the triplet
 \ $(B^{(b)}, C, \nu)$ \ with the given initial value \ $y_0$ \ with \ $\PP_b$ \ as its
 unique solution.
By Proposition \ref{Pro_jump_CIR}, the SDE \eqref{jump_CIR} has a pathwise unique strong
 solution (with the given deterministic initial value \ $y_0 \in \RR_+$), \ and hence
 Theorem III.2.26 in Jacod and Shiryaev \cite{JSh} yields that the set of all solutions to
 the martingale problem on the canonical space corresponding to \ $(B^{(b)}, C, \nu)$ \ has
 only one element \ $(\PP_b)$ \ yielding the desired local uniqueness.
We also mention that Theorem III.4.29 in Jacod and Shiryaev \cite{JSh} implies that under
 the probability measure \ $\PP_b$, \ all local martingales have the integral representation
 property relative to \ $\eta$.

By Theorem III.5.34 in Jacod and Shiryaev \cite{JSh} (see also Appendix \ref{App_LR}),
 \ $\PP_{b,T}$ \ and \ $\PP_{\tb,T}$ \ are equivalent (one can change the roles of \ $b$
 \ and \ $\tb$), \ and under the probability measure \ $\PP_{\tb}$, \  we have
 \[
   \frac{\dd \PP_{b,T}}{\dd \PP_{\tb,T}}(\eta)
  = \exp\bigg\{\int_0^T \beta^{(\tb,b)}_u \, \dd (\eta^\cont)^{(\tb)}_u
               -\frac{1}{2}
                \int_0^T \bigl(\beta^{(\tb,b)}_u\bigr)^2 c_u \, \dd u\bigg\} ,
   \qquad T \in \RR_{++} ,
 \]
 where \ $((\eta^\cont)^{(\tb)}_t)_{t\in\RR_+}$ \ denotes the continuous (local) martingale
 part of \ $(\eta_t)_{t\in\RR_+}$ \ under \ $\PP_{\tb}$.
\ Using part 1) of Remarks III.2.28 in Jacod and Shiryaev \cite{JSh} and
 \eqref{Grigelionis}, the continuous (local) martingale part
 \ $(\tY^\cont_t)_{t\in\RR_+}$ \ of \ $(\tY_t)_{t\in\RR_+}$ \ takes the form
 \ $\tY^\cont_t = \sigma \int_0^t \sqrt{\tY_u} \, \dd W_u$, \ $t \in \RR_+$, \ and, by
 \eqref{jump_CIR}, we have
 \[
   \dd \tY^\cont_t = \dd \tY_t - (a - \tb \tY_t) \, \dd t - \dd J_t , \qquad
   t \in \RR_+ .
 \]
Hence, under \ $\PP$,
 \begin{align*}
  \log\biggl(\frac{\dd \PP_{b,T}}{\dd \PP_{\tb,T}}(\tY)\bigg)
  &= \int_0^T \!\biggl(-\frac{b-\tb}{\sigma^2}\biggr) (\dd \tY_u - \dd J_u)
     - \int_0^T \!\biggl(-\frac{b-\tb}{\sigma^2}\biggr) (a - \tb \tY_u) \, \dd u
     - \frac{1}{2}
       \int_0^T \!\biggl(-\frac{b-\tb}{\sigma^2}\biggr)^2 \sigma^2 \tY_u \, \dd u \\
  &= - \frac{b-\tb}{\sigma^2} \int_0^T (\dd \tY_u - \dd J_u)
     + \frac{b-\tb}{\sigma^2} \int_0^T a \, \dd u
     - \frac{b^2-\tb^2}{2\sigma^2} \int_0^T \tY_u \, \dd u ,
 \end{align*}
 which yields the statement.
\proofend

Next, using Proposition \ref{RNb}, by considering \ $\PP_{\tb,T}$ \ as a fixed reference measure,
 we derive an MLE for the parameter \ $b$ \ based on the observations \ $(Y_t)_{t\in[0,T]}$.
\ Our method for deriving an MLE is one of the known ones in the literature, and it
 turns out that these lead to the same estimator \ $\hb_T$, \ see Remark \ref{Luschgy}.
Let us denote the right hand side of \eqref{RNformulab} by
 \ $\Lambda_T(b,\tb)$ \ replacing \ $\tY$ \ by \ $Y$.
\ By an MLE \ $\hb_T$ \ of the parameter \ $b$ \ based on the observations \ $(Y_t)_{t\in[0,T]}$,
 \ we mean
 \[
   \hb_T := \argmax_{b\in\RR} \Lambda_T(b,\tb) ,
 \]
 which will turn out to be not dependent on \ $\tb$.
\ Next, we formulate a lemma about the unique existence of MLE \ $\hb_T$ \ of \ $b$ \ for all
 \ $T \in \RR_{++}$.

\begin{Pro}\label{LEMMA_MLEb_exist}
Let \ $a \in \RR_+$, \ $b \in \RR$, \ $\sigma \in \RR_{++}$, \ $y_0 \in \RR_+$, \ and let
 \ $m$ \ be a L\'evy measure on \ $\RR_{++}$ \ satisfying \eqref{help_Levy}.
If \ $a \in \RR_{++}$ \ or \ $y_0 \in \RR_{++}$, \ then for each \ $T \in \RR_{++}$,
 \ there exists a unique MLE \ $\hb_T$ \ of \ $b$ \ almost surely having the form
 \begin{equation}\label{MLEb}
  \hb_T
  = - \frac{Y_T - y_0 - a T - J_T }{\int_0^T Y_s \, \dd s},
 \end{equation}
 provided that \ $\int_0^T Y_s \, \dd s \in\RR_{++}$ \ (which holds almost surely due to
 Proposition 2.1).
\end{Pro}

\noindent{\bf Proof.}
Due to Proposition \ref{Pro_jump_CIR},
 \ $\PP\bigl(\int_0^T Y_s \, \dd s  \in\RR_{++} \bigr) = 1$ \ for all \ $T \in \RR_{++}$,
 and hence the right hand side of \eqref{MLEb} is well-defined almost surely.
The aim of the following discussion is to show that the right hand side of \eqref{MLEb}
 is a measurable function of \ $(Y_u)_{u\in[0,T]}$ \ (i.e., a statistic).
By the L\'evy--It\^o's representation \eqref{J_Levy_Ito} of \ $J_t$, \ $t \in \RR_+$,
 \ we obtain
 \begin{equation}\label{reprJ}
  J_t = \sum_{s\in[0,t]} \Delta J_s , \qquad t\in \RR_+ ,
 \end{equation}
 since \ $\int_{-1}^1 |z| \, m(\dd z) = \int_0^1 z \, m(\dd z) < \infty$, \ see, e.g.,
 Sato \cite[Theorem 19.3]{Sat}.
Using the SDE \eqref{jump_CIR}, we have \ $\Delta J_t = \Delta Y_t$, \ $t \in \RR_+$, \ and
 then, by \eqref{reprJ}, we obtain
 \[
   J_t = \sum_{s\in[0,t]} \Delta Y_s , \qquad t \in \RR_+ .
 \]
Consequently, for all \ $t \in [0, T]$, \ $J_t$ \ is a measurable function of
 \ $(Y_u)_{u\in[0,T]}$, \ yielding that the right hand side of \eqref{MLEb} is a measurable
 function (i.e., a statistic) of \ $(Y_u)_{u\in[0,T]}$, \ as desired.
By Proposition \ref{RNb}, for all \ $b,\,\tb\in\RR$, \ we have
 \begin{align*}
  &\frac{\partial}{\partial b}
   \Lambda_T(b,\tb)
   = - \frac{1}{\sigma^2} (Y_T - y_0 - a T - J_T)
     - \frac{b}{\sigma^2} \int_0^T Y_s \, \dd s , \\
  &\frac{\partial^2}{\partial b^2}
   \Lambda_T(b,\tb)
   = - \frac{1}{\sigma^2} \int_0^T Y_s \, \dd s .
 \end{align*}
Thus the MLE \ $\hb_T$ \ of \ $b$ \ based on a sample \ $(Y_s)_{s\in[0,T]}$ \ exists almost
 surely, and it takes the form \eqref{MLEb} provided that
 \ $\int_0^T Y_s \, \dd s \in\RR_{++}$.
\proofend

In fact, it turned out that for the calculation of the MLE of \ $b$, \ one does not need to
 know the value of the parameter \ $\sigma \in \RR_{++}$ \ or the measure \ $m$.

\begin{Rem}\label{Luschgy}
In the literature there is another way of deriving an MLE.
S{\o}rensen \cite{SorM} defined an MLE of \ $\bpsi$ \ as a solution of
 the equation \ $\dot{\Lambda}_T(\bpsi) = 0$, \ where \ $\dot{\Lambda}_T(\bpsi)$ \ is the so-called  score vector
 given in formula (3.3) in S{\o}rensen \cite{SorM}.
Luschgy \cite{Lus2}, \cite{Lus} called this equation as an estimating equation.
With the notations of the proof of Proposition \ref{RNb}, taking into account of the form of
 \ $\beta^{(\tb,b)}$ \ and the fact that \ $V^{(\tb,b)}$ \ is identically one, we have
 \begin{align*}
  \dot{\Lambda}_T(b)
  &:= \int_0^T \biggl(-\frac{1}{\sigma^2}\biggr) \, \dd Y^\cont_u
   = -\frac{1}{\sigma^2} \int_0^T (\dd Y_u - (a - b Y_u) \, \dd u - \dd J_u) \\
  &= -\frac{1}{\sigma^2} \biggl(Y_T - y_0 - a T + b \int_0^T Y_u \, \dd u - J_T\biggr)
 \end{align*}
 for \ $b \in \RR$ \ and \ $T \in (0, \infty)$.
\ The estimating equation \ $\dot{\Lambda}_T(b) = 0$, \ $b \in \RR$, \ has a unique solution
 \ $-\frac{Y_T-y_0-aT-J_T}{\int_0^T Y_u\,\dd u}$ \ provided that \ $\int_0^T Y_u\,\dd u$ \ is strictly
 positive, which holds almost surely.
Recall that this unique solution coincides with \ $\hb_T$, \ see \eqref{MLEb}.
\proofend
\end{Rem}

\section{Asymptotic behaviour of the MLE in the subcritical case}
\label{section_MLE_subcritical}

\begin{Thm}\label{Thm_Laplace_subcritical}
Let \ $a \in \RR_+$, \ $b \in \RR_{++}$, \ $\sigma \in \RR_{++}$, \ and let \ $m$ \ be a
 L\'evy measure on \ $\RR_{++}$ \ satisfying \eqref{help_Levy}.
Let \ $(Y_t)_{t\in\RR_+}$ \ be the unique strong solution of the SDE \eqref{jump_CIR}
 satisfying \ $\PP(Y_0 = y_0) = 1$ \ with some \ $y_0 \in \RR_+$.
\ Then
 \begin{align}\label{help_limit_MLE_subcritical}
  \frac{1}{t} \int_0^t Y_s \, \dd s
  \stoch \frac{1}{b} \biggl(a + \int_0^\infty z \, m(\dd z)\biggr) \in \RR_+
  \qquad \text{as \ $t \to \infty$.}
 \end{align}
\end{Thm}

\noindent{\bf Proof.}
Using Corollary \ref{Thm_Laplace_intY}, we have
 \[
   \EE\left[\exp\left\{\frac{v}{t} \int_0^t Y_s \, \dd s\right\}\right]
   = \exp\bigg\{\psi_{0,\frac{v}{t}}(t) y_0
                + \int_0^t
                   \biggl(a \psi_{0,\frac{v}{t}}(s)
                          + \int_0^\infty
                             \bigl(\ee^{z\psi_{0,\frac{v}{t}}(s)}
                                   - 1\bigr) m(\dd z)\biggr) \dd s \bigg\}
 \]
 for \ $t \in \RR_{++}$ \ and \ $v \in \RR_-$, \ where the function
 \ $\psi_{0,v} : \RR_+ \to \RR_-$ \ is given in \eqref{psi_form_intY}.
We check that
 \[
  \EE\left[\exp\left\{\frac{v}{t} \int_0^t Y_s \, \dd s\right\}\right]
     \to \exp\left\{ \frac{v}{b} \left(a + \int_0^\infty z\,m(\dd z)\right)\right\}
    \qquad \text{as \ $t\to\infty$}
 \]
 for all \ $v\in\RR_{-}$.
\ The case of \ $v = 0$ \ is trivial, so we may and do suppose that \ $v \in \RR_{--}$.
\ Then we have
 \begin{align*}
  \psi_{0,\frac{v}{t}}(t)
  = \frac{1}{t}
    \frac{2v(1-\ee^{-t\gamma_{\frac{v}{t}}})}
         {\gamma_{\frac{v}{t}}(1+\ee^{-t\gamma_{\frac{v}{t}}})
          +b(1-\ee^{-t\gamma_{\frac{v}{t}}})}
  \to 0 \qquad \text{as \ $t \to \infty$,}
 \end{align*}
 since \ $\gamma_{\frac{v}{t}} = \sqrt{b^2 - \frac{2\sigma^2v}{t}} \to b \in \RR_{++}$ \ as
 \ $t \to \infty$.
\ Moreover, by the dominated convergence theorem,
 \begin{align*}
  \int_0^t \psi_{0,\frac{v}{t}}(s) \, \dd s
  &= \frac{1}{t}
     \int_0^t
      \frac{2v(1-\ee^{-s\gamma_{\frac{v}{t}}})}
           {\gamma_{\frac{v}{t}}(1+\ee^{-s\gamma_{\frac{v}{t}}})
            +b(1-\ee^{-s\gamma_{\frac{v}{t}}})}
      \dd s \\
  &= \int_0^1
      \frac{2v(1-\ee^{-xt\gamma_{\frac{v}{t}}})}
           {\gamma_{\frac{v}{t}}(1+\ee^{-xt\gamma_{\frac{v}{t}}})
            +b(1-\ee^{-xt\gamma_{\frac{v}{t}}})}
      \, \dd x
   \to \frac{v}{b} \qquad \text{as \ $t \to \infty$.}
 \end{align*}
Indeed, the integrand in the last integral is dominated by \ $\frac{2|v|}{b}$.
\ Finally, by the dominated convergence theorem, we obtain
 \begin{align*}
  &\int_0^t
   \biggl(\int_0^\infty
           \bigl(\ee^{z\psi_{0,\frac{v}{t}}(s)} - 1\bigr)
           m(\dd z)\biggr) \dd s \\
  &= \int_0^t
      \biggl(\int_0^\infty
              \biggl(\exp\biggl\{\frac{z}{t}\cdot
                                 \frac{2v(1-\ee^{-s\gamma_{\frac{v}{t}}})}
                                      {\gamma_{\frac{v}{t}}(1+\ee^{-s\gamma_{\frac{v}{t}}})
                                       +b(1-\ee^{-s\gamma_{\frac{v}{t}}})}\biggr\}
                     - 1\biggr)
              m(\dd z)\biggr) \dd s \\
  &= \int_0^1
      \biggl(\int_0^\infty
              t\biggl(\exp\biggl\{\frac{z}{t}\cdot
                                  \frac{2v(1-\ee^{-xt\gamma_{\frac{v}{t}}})}
                                      {\gamma_{\frac{v}{t}}(1+\ee^{-xt\gamma_{\frac{v}{t}}})
                                        +b(1-\ee^{-xt\gamma_{\frac{v}{t}}})}
                                 \biggr\}
                     - 1\biggr)
              m(\dd z)\biggr) \dd x\\
  &\to \int_0^1 \biggl(\int_0^\infty z \frac{v}{b} \, m(\dd z)\biggr) \dd x
   = \frac{v}{b} \int_0^\infty z \, m(\dd z)
 \end{align*}
 as \ $t \to \infty$ \ for all \ $z \in \RR_+$.
\ Indeed,
 \[
   A_t(z, x)
   := z\cdot
      \frac{2v(1-\ee^{-xt\gamma_{\frac{v}{t}}})}
           {\gamma_{\frac{v}{t}}(1+\ee^{-xt\gamma_{\frac{v}{t}}})
            +b(1-\ee^{-xt\gamma_{\frac{v}{t}}})}
   \to z \frac{v}{b} \qquad \text{as \ $t \to \infty$,}
 \]
 where \ $A_t(z, x) \in \RR_-$ \ for all \ $t \in \RR_{++}$, \ $z \in \RR_+$ \ and
 \ $x \in [0, 1]$.
\ Hence \ $t (\ee^{\frac{A_t(z, x)}{t}} - 1) \to z \frac{v}{b}$ \ as \ $t \to \infty$ \ for
 all \ $z \in \RR_+$ \ and \ $x \in [0, 1]$, \ and
 \ $|t (\ee^{\frac{A_t(z, x)}{t}} - 1)| \leq |A_t(z, x)| \leq \frac{2z|v|}{b}$
 \ for \ $t\in\RR_{++}$, \ $z\in\RR_+$ \ and \ $x\in[0,1]$,
 \ where \ $\frac{2z|v|}{b}$ \ is integrable with respect to the measure
 \ $m(\dd z) \, \dd x$ \ on \ $\RR_{++} \times [0, 1]$.
\ By the continuity theorem, we obtain
 \[
   \frac{1}{t} \int_0^t Y_s \, \dd s
   \distr \frac{1}{b} \biggl(a + \int_0^\infty z \, m(\dd z)\biggr) \in \RR_+
   \qquad \text{as \ $t \to \infty$,}
 \]
 which implies \eqref{help_limit_MLE_subcritical}.
\proofend

If \ $a \in \RR_{++}$ \ or \ $y_0 \in \RR_{++}$, \ then, using \eqref{MLEb} and the SDE
 \eqref{jump_CIR}, we get
 \begin{align}\label{MLEb-}
  \hb_T - b
  = - \frac{Y_T - y_0 - a T - J_T + b \int_0^T Y_s \, \dd s}{\int_0^T Y_s \, \dd s}
  = - \frac{\sigma \int_0^T \sqrt{Y_s} \, \dd W_s}{\int_0^T Y_s \, \dd s}
 \end{align}
 provided that \ $\int_0^T Y_s \, \dd s \in\RR_{++}$, \ which holds almost surely due to
 Proposition \ref{Pro_jump_CIR}.
Here note that \ $\sigma \int_0^T \sqrt{Y_s} \, \dd W_s = Y^{\mathrm{cont}}_T$,
 \ $T \in \RR_+$, \ see Remark \ref{Thm_MLE_cons_sigma}.

Despite the mistake in the formula (4.23) in Mai \cite{Mai} for the MLE \ $\hb_T$ \ of \ $b$,
\ Mai \cite[Theorem 4.3.1]{Mai} formulated the right asymptotic behavior of
 \ $\hb_T$ \ as \ $T \to \infty$, \ namely, asymptotic normality assuming
 ergodicity of the process \ $Y$.
\ We give a correct proof of this result, and in fact, we extend it as well,
 since we do not assume the ergodicity of \ $Y$.

\begin{Thm}\label{Thm_MLEb_subcritical}
Let \ $a \in \RR_{++}$, \ $b \in \RR_{++}$, \ $\sigma \in \RR_{++}$, \ and let \ $m$
 \ be a L\'evy measure on \ $\RR_{++}$ \ satisfying \eqref{help_Levy}.
Let \ $(Y_t)_{t\in\RR_+}$ \ be the unique strong solution of the SDE \eqref{jump_CIR}
 satisfying \ $\PP(Y_0 = y_0) = 1$ \ with some \ $y_0 \in \RR_+$.
\ Then the MLE \ $\hb_T$ \ of \ $b$ \ is asymptotically normal, i.e.,
 \begin{align}\label{help_conv_subcritical}
   \sqrt{T} (\hb_T - b)
   \distr \cN\left(0, \frac{\sigma^2b}{a+\int_0^\infty z\,m(\dd z)}\right)
           = \cN\left(0, \frac{\sigma^2}{\int_0^\infty y\,\pi(\dd y)}\right)
   \qquad \text{as \ $T \to \infty$,}
 \end{align}
 where \ $\pi$ \ denotes the unique stationary distribution of \ $(Y_t)_{t\in\RR_+}$ \
 (see part (i) of Theorem \ref{Ergodicity}).
Especially, \ $\hb_T$ \ is weakly consistent, i.e., \ $\hb_T \stoch b$ \ as
 \ $T \to \infty$.
\ With a random scaling,
 \[
   \frac{1}{\sigma} \biggl(\int_0^T Y_s \, \dd s\biggr)^{1/2} (\hb_T - b)
   \distr \cN(0, 1)
   \qquad \text{as \ $T \to \infty$.}
 \]
Under the additional moment condition \eqref{EXTRA}, \ $\hb_T$ \ is strongly consistent,
 i.e., \ $\hb_T \as b$ \ as \ $T \to \infty$.
\end{Thm}

\noindent{\bf Proof.}
By Proposition \ref{LEMMA_MLEb_exist}, there exists a unique MLE \ $\hb_T$ \ of \ $b$ \ for
 all \ $T \in \RR_{++}$, \ which has the form given in \eqref{MLEb}.
By (i) of Theorem \ref{Ergodicity}, \ $(Y_t)_{t\in\RR_+}$ \ has a unique stationary
 distribution \ $\pi$ \ with
 \ $\int_0^\infty y \, \pi(\dd y)
    = \bigl(a + \int_0^\infty z \, m(\dd z)\bigr) \frac{1}{b} \in \RR_{++}$.
\ By Theorem \ref{Thm_Laplace_subcritical}, \ we have
 \ $\frac{1}{T} \int_0^T Y_s \, \dd s \stoch \int_0^\infty y \, \pi(\dd y)$ \ as
 \ $T \to \infty$.
\ Hence, since the quadratic variation process of the square integrable martingale
 \ $\bigl(\int_0^t \sqrt{Y_s}\,\dd W_s\bigr)_{t\in\RR_+}$ \ takes the form
 \ $\bigl(\int_0^t Y_s \, \dd s\bigr)_{t\in\RR_+}$,
 \ by Theorem \ref{THM_Zanten} with \ $\eta:=\left(\int_0^\infty y\,\pi(\dd y)\right)^{1/2}$
  \ and Slutsky's lemma, we have
 \[
   \sqrt{T} (\hb_T - b)
   = -\sigma \frac{\frac{1}{\sqrt{T}} \int_0^T \sqrt{Y_s} \, \dd W_s}
                  {\frac{1}{T} \int_0^T Y_s \, \dd s}
   \distr -\sigma
           \frac{\left(\int_0^\infty y\,\pi(\dd y)\right)^{1/2}\,\cN(0,1)}
                {\int_0^\infty y\,\pi(\dd y)}
   = \cN\biggl(0, \frac{\sigma^2}{\int_0^\infty y\,\pi(\dd y)}\biggr)
 \]
 as \ $T \to \infty$, \ hence we obtain \eqref{help_conv_subcritical}.
Further, Slutsky's lemma yields
 \begin{align*}
  \frac{1}{\sigma} \left(\int_0^T Y_s \, \dd s \right)^{1/2} (\hb_T - b)
  &= \frac{1}{\sigma} \left(\frac{1}{T} \int_0^T Y_s \, \dd s \right)^{1/2}
     \sqrt{T} (\hb_T - b) \\
  &\distr \frac{1}{\sigma} \left(\int_0^\infty y\,\pi(\dd y)\right)^{1/2} \,
          \cN\biggl(0, \frac{\sigma^2}{\int_0^\infty y\,\pi(\dd y)}\biggr)
   = \cN(0, 1)
 \end{align*}
 as \ $T \to \infty$.

Under the additional moment condition \eqref{EXTRA}, by (ii) of Theorem \ref{Ergodicity}, we
 have \ $\frac{1}{T} \int_0^T Y_s \, \dd s \as \int_0^\infty y \, \pi(\dd y)$ \ as
 \ $T \to \infty$, \ implying also \ $\int_0^T Y_s \, \dd s \as \infty$ \ as
 \ $T \to \infty$.
\ Using \eqref{MLEb-} and Theorem \ref{DDS_stoch_int}, we have \ $\hb_T \as b$ \ as \ $T \to \infty$.
\proofend

We note that the moment condition \eqref{EXTRA} does not imply the moment condition \eqref{help_Levy}
 in general, and, since the moment condition \eqref{help_Levy} is already needed for us for the existence
 of a pathwise unique strong solution of the SDE \eqref{jump_CIR} (see Proposition \ref{Pro_jump_CIR}),
 in order to get our strong consistency result in Theorem \ref{Thm_MLEb_subcritical} we need to assume
 both \eqref{help_Levy} and \eqref{EXTRA}.

\section{Asymptotic behaviour of the MLE in the critical case}
\label{section_MLE_critical}

\begin{Thm}\label{Thm_Laplace_critical}
Let \ $a \in \RR_+$, \ $b = 0$, \ $\sigma \in \RR_{++}$, \ and let \ $m$ \ be a
 L\'evy measure on \ $\RR_{++}$ \ satisfying \eqref{help_Levy}.
Let \ $(Y_t)_{t\in\RR_+}$ \ be the unique strong solution of the SDE \eqref{jump_CIR}
 satisfying \ $\PP(Y_0 = y_0) = 1$ \ with some \ $y_0 \in \RR_+$.
\ Then
 \begin{align}\label{help_limit_MLE_critical2}
  \left(\frac{Y_t}{t}, \frac{1}{t^2} \int_0^t Y_s \, \dd s\right)
  \distr \left(\cY_1,\int_0^1\cY_s\,\dd s\right)
  \qquad \text{as \ $t \to \infty$,}
 \end{align}
 where \ $(\cY_t)_{t\in\RR_+}$ \ is the unique strong solution of a critical (diffusion type)
 CIR model
 \begin{align}\label{help_limit_MLE_critical}
  \dd \cY_t
  = \left(a + \int_0^\infty z \, m(\dd z)\right) \, \dd t
    + \sigma \sqrt{\cY_t} \, \dd \cW_t ,
  \qquad t \in\RR_+,
 \end{align}
 with initial condition \ $\cY_0 = 0$, \ where \ $(\cW_t)_{t\in\RR_+}$ \ is a 1-dimensional
 standard Wiener process.
Moreover, the Laplace transform of \ $(\cY_1, \int_0^1 \cY_s \, \dd s)$ \ takes the form
 \begin{align}\label{help_Laplace_critical}
  \EE\bigl(\ee^{u\cY_1+v\int_0^1\cY_s\,\dd s}\bigr)
  = \begin{cases}
     \Bigl(\cosh\left(\frac{\gamma_v}{2}\right)
           - \frac{\sigma^2 u}{\gamma_v}
             \sinh\left(\frac{\gamma_v}{2}\right)\Bigr)^{-\frac{2}{\sigma^2}
                                                          \left(a+\int_0^\infty
                                                                   z\,m(\dd z)\right)}
     & \text{if \ $v \in \RR_{--}$,} \\[2mm]
     \left(1 - \frac{\sigma^2u}{2}\right)^{-\frac{2}{\sigma^2}
                                            \left(a+\int_0^\infty z\,m(\dd z)\right)}
     & \text{if \ $v = 0$,}\\
  \end{cases}
 \end{align}
 for all \ $u, v \in \RR_-$, \ where \ $\gamma_v = \sqrt{- 2 \sigma^2 v}$, \ $v \in \RR_-$ \ (since now \ $b=0$).
\end{Thm}

\begin{Rem}
Under the conditions of Theorem \ref{Thm_Laplace_critical}, in the special case of
 \ $m = 0$, \ i.e., in the diffusion case, we get back part 1 of Theorem 1 in Ben Alaya and
 Kebaier \cite{BenKeb2}.
\proofend
\end{Rem}

\noindent{\bf Proof of Theorem \ref{Thm_Laplace_critical}.}
Using Theorem \ref{Thm_Laplace_joint}, we have
 \[
   \EE\left[\exp\left\{u \frac{Y_t}{t} + v \frac{1}{t^2} \int_0^t Y_s \, \dd s\right\}\right]
   = \exp\bigg\{\psi_{\frac{u}{t},\frac{v}{t^2}}(t) y_0
                + \int_0^t
                   \biggl(a \psi_{\frac{u}{t},\frac{v}{t^2}}(s)
                          + \int_0^\infty
                             \bigl(\ee^{z\psi_{\frac{u}{t},\frac{v}{t^2}}(s)}
                                   - 1\bigr) m(\dd z)\biggr) \dd s \bigg\}
 \]
 for \ $t \in \RR_+$ \ and \ $u, v \in \RR_-$, \ where the function
 \ $\psi_{u,v} : \RR_+ \to \RR_-$ \ is given in \eqref{psi_form}.
According to Theorem \ref{Thm_Laplace_joint}, we have to distinguish two cases:
 \ $v < 0$ \ (Case I) \ and \ $v = 0$ \ (Case II).

Case I: If \ $v < 0$, \ then
 \begin{align*}
  \psi_{\frac{u}{t},\frac{v}{t^2}}(t)
  = \frac{1}{t}
    \frac{u\gamma_v \cosh\left(\frac{\gamma_v}{2}\right)
          +2v\sinh\left(\frac{\gamma_v}{2}\right)}
         {\gamma_v\cosh\left(\frac{\gamma_v}{2}\right)
          -\sigma^2 u\sinh\left(\frac{\gamma_v}{2}\right)}
  \to 0 \qquad \text{as \ $t \to \infty$,}
 \end{align*}
 and, for each \ $t \in \RR_+$,
 \begin{align}\label{cI}
  \begin{aligned}
   \int_0^t \psi_{\frac{u}{t},\frac{v}{t^2}}(s) \, \dd s
   &= \frac{1}{t}
      \int_0^t
       \frac{u \gamma_v \cosh\left(\frac{\gamma_vs}{2t}\right)
             + 2 v \sinh\left(\frac{\gamma_vs}{2t}\right)}
            {\gamma_v \cosh\left(\frac{\gamma_vs}{2t}\right)
             - \sigma^2 u \sinh\left(\frac{\gamma_vs}{2t}\right)}
       \dd s \\
   &= \int_0^1
       \frac{u\gamma_v \cosh\left(\frac{\gamma_vx}{2}\right)
             +2v\sinh\left(\frac{\gamma_vx}{2}\right)}
            {\gamma_v \cosh\left(\frac{\gamma_vx}{2}\right)
             - \sigma^2 u \sinh\left(\frac{\gamma_vx}{2}\right)}
       \, \dd x \\
   &= \int_0^1 \psi_{u,v}(x) \, \dd x
    = -\frac{2}{\sigma^2}
       \log\Bigl(\cosh\left(\frac{\gamma_v}{2}\right)
                 -\frac{\sigma^2 u}{\gamma_v}\sinh\left(\frac{\gamma_v}{2}\right)\Bigr)
    =: \cI ,
  \end{aligned}
 \end{align}
 where we used \eqref{intpsi}.
By the monotone convergence theorem we obtain
 \begin{align*}
  &\int_0^t
   \biggl(\int_0^\infty
           \bigl(\ee^{z\psi_{\frac{u}{t},\frac{v}{t^2}}(s)} - 1\bigr)
           m(\dd z)\biggr) \dd s \\
  &= \int_0^t
      \biggl(\int_0^\infty
              \biggl(\exp\biggl\{\frac{z}{t}\cdot
                                 \frac{u \gamma_v \cosh\left(\frac{\gamma_vs}{2t}\right)
           + 2 v \sinh\left(\frac{\gamma_vs}{2t}\right)}
          {\gamma_v \cosh\left(\frac{\gamma_vs}{2t}\right)
           - \sigma^2 u \sinh\left(\frac{\gamma_vs}{2t}\right)}\biggr\}
                     - 1\biggr)
              m(\dd z)\biggr) \dd s \\
  &= \int_0^1
      \biggl(\int_0^\infty
              t\biggl(\exp\biggl\{\frac{z}{t}\cdot
                                 \frac{u \gamma_v \cosh\left(\frac{\gamma_vx}{2}\right)
           + 2 v \sinh\left(\frac{\gamma_vx}{2}\right)}
          {\gamma_v \cosh\left(\frac{\gamma_vx}{2}\right)
           - \sigma^2 u \sinh\left(\frac{\gamma_vx}{2}\right)}
                                  \biggr\}
                     - 1\biggr)
              m(\dd z)\biggr) \dd x\\
  &= \int_0^1
      \biggl(\int_0^\infty
              t(e^{\frac{z}{t}\psi_{u,v}(x)} - 1) \,
              m(\dd z)\biggr) \dd x
   \to \cI \int_0^\infty z \, m(\dd z)
 \end{align*}
 as \ $t \to \infty$ \ for all \ $z \in \RR_+$, \ where we used that
 \ $\psi_{u,v}(t) \in \RR_-$ \ for all \ $t \in \RR_+$ \ and \ $u, v \in \RR_-$, \ and
 the function \ $\RR_{++} \ni t \mapsto t (\ee^{\frac{A}{t}} - 1)$ \ is monotone decreasing
 for all \ $A \in \RR_{-}$.
\ The equality \eqref{help_Laplace_critical} is a consequence of Theorem
 \ref{Thm_Laplace_joint} and \eqref{cI}, and hence continuity theorem yields
 \eqref{help_limit_MLE_critical2}.

Case II: If \ $v =0 $, \ then
 \begin{align*}
  \psi_{\frac{u}{t},\frac{v}{t^2}}(t)
  = \frac{1}{t} \frac{u}{1-\frac{\sigma^2 u}{2}}
  \to 0 \qquad \text{as \ $t \to \infty$,}
 \end{align*}
 and, for each \ $t \in \RR_+$,
 \begin{align}\label{cII}
  \int_0^t \psi_{\frac{u}{t},\frac{v}{t^2}}(s) \, \dd s
  = \frac{1}{t}
    \int_0^t \frac{u}{1-\frac{\sigma^2 us}{2t}} \, \dd s
  = \int_0^1 \frac{u}{1-\frac{\sigma^2 ux}{2}} \, \dd x
  = - \frac{2}{\sigma^2} \log\left(1 - \frac{\sigma^2u}{2}\right) .
 \end{align}
Moreover, by the monotone convergence theorem,
 \begin{align*}
  \int_0^t
   \biggl(\int_0^\infty
           \bigl(\ee^{z\psi_{\frac{u}{t},\frac{v}{t^2}}(s)} - 1\bigr) m(\dd z)\biggr) \dd s
  &= \int_0^t
      \biggl(\int_0^\infty
              \biggl(\exp\biggl\{\frac{z}{t}\cdot
                                 \frac{u}{1-\frac{\sigma^2 us}{2t}}\biggr\}
                     - 1\biggr)
              m(\dd z)\biggr) \dd s \\
  &= \int_0^1
      \biggl(\int_0^\infty
              t\biggl(\exp\biggl\{\frac{z}{t}\cdot
                                  \frac{u}{1-\frac{\sigma^2 ux}{2}}
                                  \biggr\}
                     - 1\biggr)
              m(\dd z)\biggr) \dd x \\
  &\to - \frac{2}{\sigma^2} \log\left(1 - \frac{\sigma^2u}{2}\right)
         \int_0^\infty z \, m(\dd z) \qquad \text{as \ $t \to \infty$.}
 \end{align*}
The equality \eqref{help_Laplace_critical} is a consequence of Theorem
 \ref{Thm_Laplace_joint} and \eqref{cII}, and hence continuity theorem yields
 \eqref{help_limit_MLE_critical2}.
\proofend

\begin{Thm}\label{Thm_MLE_critical_spec}
Let \ $a \in \RR_+$, \ $b = 0$, \ $\sigma \in \RR_{++}$, \ and let \ $m$ \ be a
 L\'evy measure on \ $\RR_{++}$ \ satisfying \eqref{help_Levy}.
Let \ $(Y_t)_{t\in\RR_+}$ \ be the unique strong solution of the SDE \eqref{jump_CIR}
 satisfying \ $\PP(Y_0 = y_0) = 1$ \ with some \ $y_0 \in \RR_+$.
\ Suppose that \ $a \in \RR_{++}$ \ or \ $a=0$, \ $y_0 \in \RR_{++}$, \ $\int_0^\infty z\,m(\dd z)\in\RR_{++}$.
\ Then
 \begin{align}\label{aabb}
  T (\hb_T - b) = T \hb_T
  \distr \frac{a + \int_0^\infty  z \, m(\dd z) - \cY_1}{\int_0^1 \cY_s \, \dd s}
  \qquad \text{as \ $T \to \infty$,}
 \end{align}
 where \ $(\cY_t)_{t\in\RR_+}$ \ is the unique strong solution of the SDE \eqref{help_limit_MLE_critical}.
As a consequence, the MLE \ $\hb_T$ \ of \ $b$ \ is weakly consistent, i.e.,
 \ $\hb_T$ \ converges to \ $b$ \ in probability as \ $T\to\infty$.

With a random scaling, we have
 \[
   \frac{1}{\sigma} \biggl(\int_0^T Y_s \, \dd s\biggr)^{1/2} (\hb_T - b)
   = \frac{1}{\sigma} \biggl(\int_0^T Y_s \, \dd s\biggr)^{1/2} \, \hb_T
   \distr
   \frac{a + \int_0^\infty z \, m(\dd z) - \cY_1}
        {\sigma\bigl(\int_0^1 \cY_s \, \dd s\bigr)^{1/2}}
   \qquad \text{as \ $T \to \infty$.}
 \]
\end{Thm}

\noindent{\bf Proof.}
By Proposition \ref{LEMMA_MLEb_exist}, there exists a unique MLE \ $\hb_T$ \ of \ $b$ \ for
 all \ $T \in \RR_{++}$, \ which has the form given in \eqref{MLEb}.
By \eqref{MLEb}, we have
 \begin{align*}
  T \hb_T
  = - \frac{\frac{Y_T}{T} - \frac{y_0}{T} - a - \frac{J_T}{T}}
           {\frac{1}{T^2} \int_0^T Y_s \, \dd s} ,
  \qquad T \in\RR_+ .
 \end{align*}
Here, by strong law of large numbers for the L\'evy process \ $(J_t)_{t\in\RR_+}$ \ (see,
 e.g., Kyprianou \cite[Exercise 7.2]{Kyp}),
 \begin{align}\label{help_Levy_SLLN}
  \PP\left(\lim_{T\to\infty} \frac{J_T}{T} = \EE(J_1) = \int_0^\infty z \, m(\dd z)\right)
  = 1 .
 \end{align}
Hence, using Slutsky's lemma, Theorem \ref{Thm_Laplace_critical} and the continuous mapping
 theorem together with the fact that \ $\PP(\int_0^1 \cY_s\,\dd s\in\RR_{++})=1$ \
 (following from Proposition \ref{Pro_jump_CIR}, since under the conditions of
 Theorem \ref{Thm_MLE_critical_spec}, \ $a + \int_0^\infty z \, m(\dd z) \in\RR_{++}$),
 \ we have
 \[
   T \hb_T \distr -\frac{\cY_1 - a - \int_0^\infty z \, m(\dd z)}{\int_0^1 \cY_s \, \dd s}
   \qquad \text{as \ $T \to \infty$,}
 \]
 as desired.
Consequently, by Slutsky's lemma, weak consistency follows by the decomposition
 \ $\hb_T = \frac{1}{T} T \hb_T$, \ $T \in \RR_{++}$.

Applying Slutsky's lemma, Theorem \ref{Thm_Laplace_critical} and the continuous mapping
 theorem, we obtain
 \begin{align*}
  &\frac{1}{\sigma} \left(\int_0^T Y_s \, \dd s\right)^{1/2} \, \hb_T
   = - \frac{1}{\sigma} \left(\frac{1}{T^2} \int_0^T Y_s \, \dd s \right)^{-1/2} \,
       \left(\frac{Y_T}{T} - \frac{y_0}{T} - a - \frac{J_T}{T}\right) \\
  &\distr - \frac{1}{\sigma} \left(\int_0^1 \cY_s \, \dd s\right)^{-1/2} \,
           \left(\cY_1 - a - \int_0^\infty z \, m(\dd z)\right)
   =\frac{a + \int_0^\infty z \, m(\dd z) - \cY_1}
         {\sigma\bigl(\int_0^1 \cY_s \, \dd s\bigr)^{1/2}}
   \qquad \text{as \ $T \to \infty$,}
 \end{align*}
 as desired.
\proofend

\section{Asymptotic behaviour of the MLE in the supercritical case}
\label{section_MLE_supercritical}

\begin{Thm}\label{Thm_supercritical_convergence}
Let \ $a \in \RR_+$, \ $b \in \RR_{--}$, \ $\sigma \in \RR_{++}$, \ and let \ $m$ \ be a
 L\'evy measure on \ $\RR_{++}$ \ satisfying \eqref{help_Levy}.
Let \ $(Y_t)_{t\in\RR_+}$ \ be the unique strong solution of the SDE \eqref{jump_CIR}
 satisfying \ $\PP(Y_0 = y_0) = 1$ \ with some \ $y_0 \in \RR_+$.
\ Then there exists a random variable \ $V$ \ with \ $\PP(V \in \RR_+) = 1$ \ such that
 \[
   \ee^{bt} Y_t \as V \qquad \text{and} \qquad
   \ee^{bt} \int_0^t Y_u \, \dd u \as -\frac{V}{b}
   \qquad \text{as \ $t \to \infty$.}
 \]
Moreover, the Laplace transform of \ $V$ \ takes the form
 \begin{align}\label{Laplace_supercritical_limit}
   \EE(\ee^{u V})
      = \exp\biggl\{\frac{uy_0}{1+\frac{\sigma^2 u}{2b}}\biggr\}
        \left(1 + \frac{\sigma^2 u}{2b}\right)^{-\frac{2a}{\sigma^2}}
        \exp\biggl\{\int_0^\infty
                     \biggl(\int_0^\infty
                             \biggl(\exp\biggl\{\frac{zu\ee^{by}}
                                                     {1+\frac{\sigma^2 u}{2b}
                                                        \ee^{by}}\biggr\}
                                    - 1\biggr)
                             m(\dd z)\biggr)
                     \dd y\biggr\}
 \end{align}
 for all \ $u \in \RR_-$, \ and consequently,
 \[
   V \distre \tcV + \ttcV ,
 \]
 where \ $\tcV$ \ and \ $\ttcV$ \ are independent random variables such that
 \ $\ee^{bt} \tcY_t \as \tcV$ \ and \ $\ee^{bt} \ttcY_t \as \ttcV$ \ as
 \ $t \to \infty$, \ where \ $(\tcY_t)_{t\in\RR_+}$ \ and  \ $(\ttcY_t)_{t\in\RR_+}$
 \ are the pathwise unique strong solutions of the supercritical CIR models
 \begin{align}\label{SDE_CIR_tcY}
  \dd \tcY_t = (a - b \tcY_t) \, \dd t + \sigma \sqrt{\tcY_t} \, \dd \tcW_t ,
  \qquad t \in \RR_+ ,  \qquad \text{with}\qquad \tcY_0 = y_0 ,
 \end{align}
 and
 \begin{align}\label{SDE_CIR_ttcY}
  \dd \ttcY_t = - b \ttcY_t \, \dd t + \sigma \sqrt{\ttcY_t} \, \dd \ttcW_t + \dd J_t ,
  \qquad t \in\RR_+ ,  \qquad \text{with} \qquad \ttcY_0 = 0 ,
 \end{align}
 respectively, where \ $(\tcW_t)_{t\in\RR_+}$ \ and \ $(\ttcW_t)_{t\in\RR_+}$ \ are
 independent 1-dimensional standard Wiener processes (yielding that
 \ $(\tcY_t)_{t\in\RR_+}$ \ and \ $(\ttcY_t)_{t\in\RR_+}$ \ are independent).
Further,
 \[
   \tcV \distre \cZ_{-\frac{1}{b}} ,
 \]
 where \ $(\cZ_t)_{t\in\RR_+}$ \ is the pathwise unique strong solutions of the critical CIR
 model
 \begin{align}\label{SDE_CIR_cY}
  \dd \cZ_t = a \, \dd t + \sigma \sqrt{\cZ_t} \, \dd \cW_t,
  \qquad t \in\RR_+ ,  \qquad \text{with}\qquad \cZ_0 = y_0 ,
 \end{align}
 where \ $(\cW_t)_{t\in\RR_+}$ \ is a 1-dimensional standard Wiener process.

Moreover, \ $\PP(V \in \RR_{++}) = 1$ \ if and only if \ $a \in \RR_{++}$ \ or \ $m \ne 0$.

If, in addition, \ $a \in \RR_{++}$, \ then \ $V$ \ is absolutely continuous.
\end{Thm}

\begin{Rem}\label{super}
(i) Under the conditions of Theorem \ref{Thm_supercritical_convergence} in the special case
 of \ $m$ \ having the form \eqref{Levy_cond}, Corollary \ref{Thm_Laplace_Y} yields
 another representation of the law of \ $V$ \ in Theorem
 \ref{Thm_supercritical_convergence}, namely,
 \begin{align}\label{supercrit_stoch_rep_spec}
   V \distre \tcX_{-\frac{1}{b}}
             + \ttcX_{-\frac{1}{b}\left(1-\frac{2b}{\sigma^2 \lambda}\right)} ,
 \end{align}
 where \ $(\tcX_t)_{t\in\RR_+}$ \ and  \ $(\ttcX_t)_{t\in\RR_+}$ \ are the pathwise unique
 strong solutions of the critical CIR models
 \begin{align*}
  \dd \tcX_t = a \, \dd t + \sigma \sqrt{\tcX_t} \, \dd \tcW_t ,
  \qquad t \in \RR_+ ,  \qquad \text{with} \qquad \tcX_0 = y_0 ,
 \end{align*}
 and
 \begin{align*}
  \dd \ttcX_t = \frac{c}{\lambda-\frac{2b}{\sigma^2}} \, \dd t
                + \sigma \sqrt{\ttcX_t} \, \dd \ttcW_t ,
  \qquad t \in \RR_+ , \qquad \text{with} \qquad \ttcX_0 = 0 ,
 \end{align*}
 respectively, where \ $(\tcW_t)_{t\in\RR_+}$ \ and \ $(\ttcW_t)_{t\in\RR_+}$ \ are
 independent 1-dimensional standard Wiener processes (yielding that \ $(\tcX_t)_{t\in\RR_+}$
 \ and \ $(\ttcX_t)_{t\in\RR_+}$ \ are independent).
We point out that in the representation \eqref{supercrit_stoch_rep_spec} there are only diffusion-type (critical) CIR processes,
 while the representation of \ $V$ \ given in Theorem \ref{Thm_supercritical_convergence}
 contains the almost sure limit of the appropriately scaled jump-type (supercritical) CIR process \eqref{SDE_CIR_ttcY}.
Further, by Alfonsi \cite[Proposition 1.2.11]{Alf}, if \ $a \in \RR_{++}$, \ then
 \ $\tcX_{-\frac{1}{b}}$ \ is absolutely continuous, and, if \ $a = 0$ \ and
 \ $c \in \RR_{++}$, \ then
 \ $\ttcX_{-\frac{1}{b}\left(1-\frac{2b}{\sigma^2 \lambda}\right)}$ \ is absolutely
 continuous.
The independence of \ $\tcX_{-\frac{1}{b}}$ \ and
 \ $\ttcX_{-\frac{1}{b}\left(1-\frac{2b}{\sigma^2 \lambda}\right)}$ \ implies that \ $V$
 \ is absolutely continuous in both cases.

\noindent (ii) Under the conditions of Theorem \ref{Thm_supercritical_convergence}, in the
 special case \ $m = 0$, \ i.e., in the diffusion case, we get back part 2 of Theorem 3 in
 Ben Alaya and Kebaier \cite{BenKeb2}.
Indeed, if \ $m = 0$, \ then the pathwise unique strong solution of the SDE
 \eqref{SDE_CIR_ttcY} is the identically zero process.
\proofend
\end{Rem}

\noindent{\bf Proof of Theorem \ref{Thm_supercritical_convergence}.}
First, we prove the existence of an appropriate non-negative random variable \ $V$.
\ We check that
 \begin{align*}
  \EE(Y_t \mid \cF^Y_s)
  = \EE(Y_t \mid Y_s)
  = \ee^{-b(t-s)} Y_s
    + \biggl(a + \int_0^\infty z \, m(\dd z)\biggr) \ee^{bs} \int_s^t \ee^{-bu} \, \dd u
 \end{align*}
 for all \ $s, t \in \RR_+$ \ with \ $0 \leq s \leq t$, \ where
 \ $\cF_t^Y := \sigma(Y_s, s\in[0,t])$, \ $t \in \RR_+$.
\ The first equality follows from the Markov property of the process \ $(Y_t)_{t\in\RR_+}$.
\ The second equality is a consequence of the time-homogeneity of the Markov process \ $Y$
 \ and the fact that
 \begin{align*}
  \EE(Y_t \mid Y_0 = y_0)
  = \ee^{-bt} y_0
    + \biggl(a + \int_0^\infty z \, m(\dd z)\biggr) \int_0^t \ee^{-bu} \, \dd u ,
  \qquad t \in \RR_+ , \quad y_0\in \RR_+,
 \end{align*}
 following from  Proposition \ref{Pro_moments}.
Then
 \[
   \EE(\ee^{bt} Y_t \mid \cF^Y_s)
   = \ee^{bs} Y_s
     + \biggl(a + \int_0^\infty z \, m(\dd z)\biggr)
       \ee^{b(s+t)} \int_s^t \ee^{-bu} \, \dd u
   \geq \ee^{bs} Y_s
 \]
 for all \ $s, t \in \RR_+$ \ with \ $0 \leq s \leq t$, \ consequently, the process
 \ $(\ee^{bt} Y_t)_{t\in\RR_+}$ \ is a non-negative submartingale with respect to the
 filtration \ $(\cF_t^Y)_{t\in\RR_+}$.
\ Moreover,
 \begin{align*}
  \EE(\ee^{bt} Y_t)
  & = y_0 + \biggl(a + \int_0^\infty z \, m(\dd z)\biggr) \ee^{bt} \int_0^t \ee^{-bu} \, \dd u
    = y_0 + \biggl(a + \int_0^\infty z \, m(\dd z)\biggr)  \int_0^t \ee^{bu} \, \dd u  \\
  & \leq  y_0 + \biggl(a + \int_0^\infty z \, m(\dd z)\biggr) \int_0^\infty \ee^{bu} \, \dd u
    = y_0 - \frac{1}{b} \biggl(a + \int_0^\infty z \, m(\dd z)\biggr)
   < \infty , \qquad t \in \RR_+ ,
 \end{align*}
 hence, by the submartingale convergence theorem, there exists a non-negative random
 variable \ $V$ \ such that
 \begin{align}\label{lim_Y}
  \ee^{bt} Y_t \as V \qquad \text{as \ $t \to \infty$.}
 \end{align}
Further, if \ $\omega \in \Omega$ \ is such that \ $\ee^{bt} Y_t(\omega) \to V(\omega)$ \ as \ $t \to \infty$,
 \ then, by the integral Toeplitz lemma \ (see K\"uchler and S{\o}rensen
 \cite[Lemma B.3.2]{KucSor}), we have
 \[
   \frac{1}{\int_0^t \ee^{-bu} \, \dd u}
   \int_0^t \ee^{-bu} (\ee^{bu} Y_u(\omega)) \, \dd u
   \to V(\omega)  \qquad \text{as \ $t \to \infty$.}
 \]
Here \ $\int_0^t \ee^{-bu} \, \dd u = \frac{\ee^{-bt} - 1}{-b}$, \ $t \in \RR_+$,
 \ thus we conclude
 \begin{align}\label{lim_intY}
  \ee^{bt} \int_0^t Y_u \, \dd u
  = \frac{1-\ee^{bt}}{-b} \frac{\int_0^t Y_u \, \dd u}{\int_0^t \ee^{-bu} \, \dd u}
  \as -\frac{V}{b} \qquad
  \text{as \ $t \to \infty$.}
 \end{align}

Now we turn to calculate the Laplace transform of \ $V$.
\ Since \ $\PP(\lim_{t\to\infty} \ee^{bt} Y_t = V) = 1$, \ we also have
 \ $\ee^{bt} Y_t \distr V$ \ as \ $t \to \infty$, \ and, by continuity theorem,
 \ $\lim_{t\to\infty} \EE(\exp\{u\ee^{bt}Y_t\}) = \ee^{u V}$ \ for all \ $u \in \RR_-$.
\ By Corollary \ref{Thm_Laplace_Y},
 \[
   \EE(\exp\{u\ee^{bt}Y_t\})
   = \exp\biggl\{\psi_{u\ee^{bt},0}(t) y_0
                 + \int_0^t
                    \biggl(a \psi_{u\ee^{bt}, 0}(s)
                           + \int_0^\infty
                              \bigl(\ee^{z\psi_{u\ee^{bt}, 0}(s)}
                                    - 1\bigr) m(\dd z)\biggr) \dd s \bigg\}
 \]
 for each \ $u \in \RR_-$ \ and \ $t \in \RR_+$, \ where we have
 \begin{align*}
  \psi_{u\ee^{bt}, 0}(t)
  = \frac{2ub}{\sigma^2 u(1-\ee^{bt})+2b}
  \to \frac{u}{1+\frac{\sigma^2 u}{2b}}
  \qquad \text{as \ $t \to \infty$.}
 \end{align*}
By the monotone convergence theorem we obtain
 \begin{align*}
  \int_0^t \psi_{u\ee^{bt}, 0}(s) \, \dd s
  &= \int_0^t \frac{2ub\ee^{b(t-s)}}{\sigma^2 u(\ee^{b(t-s)}-\ee^{bt})+2b} \, \dd s
   = \int_0^t \frac{2ub\ee^{by}}{\sigma^2 u(\ee^{by}-\ee^{bt})+2b} \, \dd y \\
  &\to \int_0^\infty \frac{2ub\ee^{by}}{\sigma^2 u\ee^{by}+2b} \, \dd y
   = - \frac{2}{\sigma^2} \log\Bigl(1+\frac{\sigma^2 u}{2b}\Bigr)
  \qquad \text{as \ $t \to \infty$,}
 \end{align*}
 where we used that the function
 \ $\RR_{++} \ni t \mapsto \frac{2ub\ee^{by}}{\sigma^2 u(\ee^{by}-\ee^{bt})+2b} \in \RR_{--}$
 \ is monotone decreasing for all \ $y \in \RR_{++}$ \ and \ $u, b \in \RR_{--}$.
\ In a similar way,
 \begin{align*}
  \int_0^t
   \biggl(\int_0^\infty \bigl(\ee^{z\psi_{u\ee^{bt}, 0}(s)} - 1\bigr) m(\dd z)\biggr) \dd s
  &= \int_0^t
      \biggl(\int_0^\infty
              \biggl(\exp\biggl\{\frac{2zub\ee^{b(t-s)}}
                                      {\sigma^2 u(\ee^{b(t-s)}-\ee^{bt})+2b}\biggr\}
                     - 1\biggr)
              m(\dd z)\biggr)
      \dd s \\
  &= \int_0^t
      \biggl(\int_0^\infty
              \biggl(\exp\biggl\{\frac{2zub\ee^{by}}
                                      {\sigma^2 u(\ee^{by}-\ee^{bt})+2b}\biggr\}
                     - 1\biggr)
              m(\dd z)\biggr)
      \dd y \\
  &\to \int_0^\infty
        \biggl(\int_0^\infty
                \biggl(\exp\biggl\{\frac{2zub\ee^{by}}
                                        {\sigma^2 u\ee^{by}+2b}\biggr\} - 1\biggr)
                m(\dd z)\biggr)
        \dd y
 \end{align*}
 as \ $t \to \infty$, \ yielding \eqref{Laplace_supercritical_limit}.

Particularly, we obtain \ $\ee^{bt} \tcY_t \as \tcV$ \ and \ $\ee^{bt} \ttcY_t \as \ttcV$
 \ as \ $t \to \infty$ \ with
 \begin{align*}
  \EE(\ee^{u \tcV})
  &= \exp\biggl\{\frac{uy_0}{1+\frac{\sigma^2 u}{2b}}\biggr\}
     \left(1 + \frac{\sigma^2 u}{2b}\right)^{-\frac{2a}{\sigma^2}} , \\
  \EE(\ee^{u \ttcV})
  &= \exp\biggl\{\int_0^\infty
                  \biggl(\int_0^\infty
                          \biggl(\exp\biggl\{\frac{zu\ee^{by}}
                                                  {1+\frac{\sigma^2 u}{2b}
                                                     \ee^{by}}\biggr\}
                                 - 1\biggr)
                          m(\dd z)\biggr)
                  \dd y\biggr\}
 \end{align*}
 for all \ $u \in \RR_-$, \ and we conclude \ $V \distre \tcV + \ttcV$.
\ The equality \ $\tcV \distre \cZ_{-\frac{1}{b}}$ \ readily follows from Corollary
 \ref{Thm_Laplace_Y}, since
 \begin{align*}
  \EE(\ee^{u\cZ_t})
  = \exp\left\{\frac{uy_0}{1-\frac{\sigma^2 u}{2}t}
               + a \int_0^t \frac{u}{1-\frac{\sigma^2 u}{2} s} \, \dd s\right\}
  = \exp\left\{\frac{uy_0}{1-\frac{\sigma^2 u}{2}t}\right\}
        \left(1 - \frac{\sigma^2 u}{2}t\right)^{-\frac{2a}{\sigma^2}}
 \end{align*}
 for \ $t\in\RR_+$, \ $u\in\RR_-$, \ and hence
 \[
   \EE(\ee^{u\cZ_{-1/b}}) = \exp\left\{ \frac{uy_0}{1+\frac{\sigma^2 u}{2b} } \right\}
         \left(1 + \frac{\sigma^2 u}{2b} \right)^{-\frac{2a}{\sigma^2}} ,
   \qquad u \in \RR_- .
 \]

The monotone convergence theorem yields
 \ $\EE(\ee^{uV}) \downarrow \EE(\bbone_{\{V=0\}}) = \PP(V = 0)$ \ as \ $u\to-\infty$.
\ We have
 \[
   \exp\biggl\{\frac{uy_0}{1+\frac{\sigma^2 u}{2b}}\biggr\}
   = \exp\biggl\{\frac{y_0}{\frac{1}{u}+\frac{\sigma^2}{2b}}\biggr\}
    \,\downarrow\, \exp\biggl\{\frac{2by_0}{\sigma^2}\biggr\}
   \in \RR_{++} \qquad \text{as \ $u\to-\infty$}
 \]
 and
 \[
   \left(1 + \frac{\sigma^2 u}{2b}\right)^{-\frac{2a}{\sigma^2}}
    \downarrow \begin{cases}
        1 & \text{if \ $a = 0$,} \\
        0 & \text{if \ $a \in \RR_{++}$}
       \end{cases}
   \qquad \text{as \ $u\to-\infty$.}
 \]
Moreover, for each \ $z \in \RR_+$ \ and \ $y \in \RR_+$, \ the function
 \[
   \RR_- \ni u \mapsto
   \exp\biggl\{\frac{zu\ee^{by}}{1+\frac{\sigma^2 u}{2b} \ee^{by}}\biggr\} - 1
   = \exp\biggl\{\frac{z\ee^{by}}{\frac{1}{u}+\frac{\sigma^2}{2b} \ee^{by}}\biggr\} - 1
   \in (-1, 0]
 \]
 is monotone increasing, hence, again by the monotone convergence theorem, we obtain
 \[
   \int_{\RR_+^2}
    \biggl(\exp\biggl\{\frac{zu\ee^{by}}
                            {1+\frac{\sigma^2 u}{2b}\ee^{by}}\biggr\} - 1\biggr)
    m(\dd z) \, \dd y
    \,\downarrow\,
   \int_{\RR_+^2} \biggl(\exp\biggl\{\frac{2bz}{\sigma^2}\biggr\} - 1\biggr) m(\dd z) \, \dd y
   = \begin{cases}
      0 & \text{if \ $m = 0$,} \\
      - \infty & \text{if \ $m \ne 0$}
     \end{cases}
 \]
 as \ $u \to -\infty$.
\ This yields
 \[
   \exp\biggl\{\int_0^\infty
                \biggl(\int_0^\infty
                        \biggl(\exp\biggl\{\frac{zu\ee^{by}}
                                                {1+\frac{\sigma^2 u}{2b}\ee^{by}}\biggr\}
                               - 1\biggr)
                        m(\dd z)\biggr)
                \dd y\biggr\}
    \,\downarrow\, \begin{cases}
        1 & \text{if \ $m = 0$,} \\
        0 & \text{if \ $m \ne 0$}
       \end{cases}
   \qquad \text{as \ $u\to-\infty$.}
 \]
Summarizing, by \eqref{Laplace_supercritical_limit}, we conclude that \ $\PP(V = 0) = 0$ \ if and only if \ $a \in \RR_{++}$ \ or
 \ $m \ne 0$.
\ Consequently, \ $\PP(V \in \RR_{++}) = 1$ \ if and only if \ $a \in \RR_{++}$ \ or
 \ $m \ne 0$.

If \ $a \in \RR_{++}$, \ then the distribution of \ $\cZ_{-\frac{1}{b}}$ \ (and hence that of \ $\tcV$) \ is absolutely
 continuous (see, e.g., Alfonsi \cite[Proposition 1.2.11]{Alf}), and consequently \ $V$ \ is
 absolutely continuous.
\proofend

\begin{Thm}\label{Thm_MLE_supercritical}
Let \ $a \in \RR_+$, \ $b \in \RR_{--}$, \ $\sigma \in \RR_{++}$, \ and let \ $m$ \ be a
 L\'evy measure on \ $\RR_{++}$ \ satisfying \eqref{help_Levy}.
Let \ $(Y_t)_{t\in\RR_+}$ \ be the unique strong solution of the SDE \eqref{jump_CIR}
 satisfying \ $\PP(Y_0 = y_0) = 1$ \ with some \ $y_0 \in \RR_+$.
\ Suppose that \ $a \in \RR_{++}$ \ or \ $a = 0$, \ $m \ne 0$ \ and \ $y_0 \in \RR_{++}$.
\ Then the MLE of \ $b$ \ is strongly consistent, i.e., \ $\hb_T \as b$ \ as
 \ $T \to \infty$.
\ Further,
 \begin{align*}
  \ee^{-bT/2} (\hb_T  - b)
  \distr \sigma Z \left(-\frac{V}{b}\right)^{-1/2}
  \qquad \text{as \ $T \to \infty$,}
 \end{align*}
 where \ $V$ \ is a positive random variable having Laplace transform given in
 \eqref{Laplace_supercritical_limit}, and \ $Z$ \ is a standard normally distributed random
 variable, independent of \ $V$.

With a random scaling, we have
 \[
   \frac{1}{\sigma} \biggl(\int_0^T Y_s \, \dd s\biggr)^{1/2} (\hb_T  - b)
   \distr \cN(0, 1) \qquad \text{as \ $T \to \infty$.}
 \]
\end{Thm}

\noindent{\bf Proof.}
By Proposition \ref{LEMMA_MLEb_exist}, there exists a unique MLE \ $\hb_T$ \ of \ $b$ \ for
 all \ $T \in \RR_{++}$ \ which takes the form given in \eqref{MLEb}.
Using \eqref{MLEb}, \eqref{help_Levy_SLLN}, Theorem \ref{Thm_supercritical_convergence} and
 the fact that \ $\lim_{T\to\infty} T \ee^{bT} = 0$, \ we have
 \begin{align*}
  \hb_T
  = \frac{-\ee^{bT}Y_T+\ee^{bT}y_0+aT\ee^{bT}+T\ee^{bT}\frac{J_T}{T}}
         {\ee^{bT}\int_0^T Y_s\,\dd s}
  \as \frac{-V}{-\frac{V}{b}}
  = b
  \qquad \text{as \ $T \to \infty$,}
 \end{align*}
 as desired.
Further, by the second equality in \eqref{MLEb-},
 \begin{align*}
  \ee^{-bT/2}(\hb_T  - b)
        = - \sigma
        \frac{\ee^{bT/2}\int_0^T\sqrt{Y_s}\,\dd W_s}{\ee^{bT}\int_0^T Y_s\,\dd s} ,
  \qquad T \in \RR_{++} .
 \end{align*}
By Theorem \ref{Thm_supercritical_convergence},
 \ $\ee^{bT} \int_0^T Y_s \, \dd s \as -\frac{V}{b}$ \ as \ $T \to \infty$, \ and using
 Theorem \ref{THM_Zanten} with \ $\eta := \left(-\frac{V}{b}\right)^{1/2}$ \ and with
 \ $v := - \frac{V}{b}$, \ we have
 \[
   \left(\ee^{bT/2} \int_0^T \sqrt{Y_s} \, \dd W_s, - \frac{V}{b}\right)
   \distr \left(\left(-\frac{V}{b}\right)^{1/2} Z, - \frac{V}{b}\right)
   \qquad \text{as \ $T\to\infty$.}
 \]
Consequently,
 \begin{equation}\label{Zan}
  \left(\ee^{bT/2} \int_0^T \sqrt{Y_s} \, \dd W_s, \ee^{bT}\int_0^T Y_s\,\dd s\right)
  \distr \left(\left(-\frac{V}{b}\right)^{1/2} Z, - \frac{V}{b}\right)
  \qquad \text{as \ $T\to\infty$.}
 \end{equation}
By the continuous mapping theorem,
 \ $\ee^{-bT/2} (\hb_T - b) \distr -\sigma Z \left(\frac{-V}{b}\right)^{-1/2}$ \ as
 \ $T \to \infty$, \ yielding the first assertion.

Applying again \eqref{Zan} and the continuous mapping theorem, we obtain
 \begin{align*}
  \frac{1}{\sigma} \biggl(\int_0^T Y_s \, \dd s\biggr)^{1/2} (\hb_T  - b)
  &= - \left(\ee^{bT} \int_0^T Y_s\,\dd s\right)^{-1/2}
       \ee^{bT/2} \int_0^T \sqrt{Y_s} \, \dd W_s \\
  &\distr - \left(-\frac{V}{b}\right)^{-1/2} \left(-\frac{V}{b}\right)^{1/2} Z
   = - Z
   \distre \cN(0, 1)
  \qquad \text{as \ $T\to\infty$,}
 \end{align*}
 as desired.
\proofend

\vspace*{3mm}

\appendix

\vspace*{5mm}

\noindent{\bf\Large Appendices}

\section{A comparison theorem in the jump process}
\label{section_comparison}

Next we prove a comparison theorem for the SDE \eqref{jump_CIR} in the jump process \ $J$.

\begin{Pro}\label{Pro_comparison_jump_CIR}
Let \ $a \in \RR_+$, \ $b \in \RR$ \ and \ $\sigma \in \RR_{++}$.
\ Suppose that the L\'evy measure \ $m$ \ on \ $\RR_{++}$ \ satisfies \eqref{help_Levy}.
Let \ $\eta_0$ \ and \ $\eta'_0$ \ be random variables independent of \ $(W_t)_{t\in\RR_+}$
 \ and \ $(J_t)_{t\in\RR_+}$ \ satisfying \ $\PP(\eta_0 \in \RR_+) = 1$ \ and
 \ $\PP(\eta'_0 \in \RR_+) = 1$.
\ Let \ $(Y_t)_{t\in\RR_+}$ \ be a pathwise unique strong solution of the SDE
 \eqref{jump_CIR} such that \ $\PP(Y_0 = \eta_0) = 1$.
Let \ $(Y'_t)_{t\in\RR_+}$ \ be a pathwise unique strong solution of the SDE
 \begin{equation}\label{CIR}
  \dd Y'_t = (a - b Y'_t) \, \dd t + \sigma \sqrt{Y'_t} \, \dd W_t , \qquad t \in \RR_+ ,
 \end{equation}
 such that \ $\PP(Y'_0 = \eta'_0) = 1$.
\ Then \ $\PP(\eta_0 \geq \eta'_0) = 1$ \ implies
 \ $\PP(\text{$Y_t \geq Y_t'$ \ for all \ $t \in \RR_+$}) = 1$.
\end{Pro}

\noindent
\textbf{Proof.} \
We follow the ideas of the proof of Theorem 3.1 of Ma \cite{Ma}, which is an adaptation of
 that of Theorem 5.5 of Fu and Li \cite{FuLi}.
There is a sequence \ $\phi_k : \RR \to \RR_+$, \ $k \in \NN$, \ of twice continuously
 differentiable functions such that
 \begin{enumerate}
  \item[(i)]
   $\phi_k(z) \uparrow z^+ := \max(z,0)$ \ as \ $k \to \infty$;
  \item[(ii)]
   $\phi_k'(z) \in [0, 1]$ \ for all \ $z \in \RR_+$ \ and \ $k \in \NN$;
  \item[(iii)]
   $\phi_k'(z) = \phi_k(z) = 0$ \ whenever \ $z \in \RR_{-}$ \ and \ $k \in \NN$;
  \item[(iv)]
   $\phi_k''(x - y) (\sqrt{x} - \sqrt{y})^2 \leq 2/k$ \ for all
    \ $x, y \in \RR_+$ \ and \ $k \in \NN$.
 \end{enumerate}
For a construction of such functions, see, e.g., the proof of Theorem 3.1 of Ma \cite{Ma}.
Let \ $Z_t := Y'_t - Y_t$ \ for all \ $t \in \RR_+$.
\ By \eqref{jump_CIR}, \eqref{J_Levy_Ito} and \eqref{CIR}, we have
 \[
   Z_t = Z_0 - b \int_0^t Z_s \, \dd s
         + \sigma \int_0^t \bigl(\sqrt{Y'_s} - \sqrt{Y_s}\bigr) \, \dd W_s
         - \int_0^t \int_0^\infty z \, \mu^J(\dd s,\dd z) , \qquad t \in \RR_+ .
 \]
For each \ $j \in \NN$, \ put
 \[
   \tau_j := \inf\Bigl\{ t \in \RR_+ : \max\{Y_t , Y'_t\} \geq j \Bigr\} .
 \]
By It\^o's formula (which can be used since \ $Y$ \ and \ $Y'$ \ are adapted to the
 augmented filtration corresponding to \ $(W_t)_{t\in\RR_+}$ \ and \ $(J_t)_{t\in\RR_+}$
 \ constructed as in Section 3 of Barczy et al.\ \cite{BarLiPap1}), we obtain
 \[
   \phi_k(Z_{t\land\tau_j}) = \phi_k(Z_0) + \sum_{\ell=1}^4 I_{j,k,\ell}(t) ,
   \qquad t \in \RR_+ ,
 \]
 with
 \begin{align*}
  &I_{j,k,1}(t)
   := - b \int_0^{t\land\tau_j} \phi_k'(Z_s) Z_s \, \dd s , \\
  &I_{j,k,2}(t)
   := \sigma
      \int_0^{t\land\tau_j}
       \phi_k'(Z_s) \bigl(\sqrt{Y'_s} - \sqrt{Y_s}\bigr)
       \, \dd W_s , \\
  &I_{j,k,3}(t)
   := \frac{\sigma^2}{2}
      \int_0^{t\land\tau_j}
       \phi_k''(Z_s) \bigl(\sqrt{Y'_s} - \sqrt{Y_s}\bigr)^2
       \, \dd s ,\\
  &I_{j,k,4}(t)
   := \int_0^{t\land\tau_j} \int_0^\infty
       [\phi_k(Z_{s-} - z) - \phi_k(Z_{s-})] \, \mu^J(\dd s, \dd z) .
 \end{align*}
Using formula (3.8) in Chapter II in Ikeda and Watanabe \cite{IkeWat}, the last integral can
 be written as \ $I_{j,k,4}(t) = I_{j,k,4,1}(t) + I_{j,k,4,2}(t)$, \ where
 \begin{align*}
  I_{j,k,4,1}(t)
  &:= \int_0^{t\land\tau_j} \int_0^\infty
        [\phi_k(Z_{s-} - z) - \phi_k(Z_{s-})] \, \tmu^J(\dd s, \dd z) , \\
  I_{j,k,4,2}(t)
  &:= \int_0^{t\land\tau_j} \int_0^\infty
       [\phi_k(Z_{s-} - z) - \phi_k(Z_{s-})] \, \dd s \, m(\dd z) ,
 \end{align*}
 where \ $\tmu^J(\dd s, \dd z) := \mu^J(\dd s, \dd z) - \dd s\,m(\dd z)$, \ since
 \begin{equation}\label{Fp1}
  \begin{aligned}
   &\EE\left( \int_0^{t\land\tau_j} \int_0^\infty
               |\phi_k(Z_{s-} - z) - \phi_k(Z_{s-})|
               \, \dd s \, m(\dd z) \right) \\
   &\leq\EE\left( \int_0^{t\land\tau_j} \int_0^\infty
                   z
                   \, \dd s \, m(\dd z) \right)
    \leq t \int_0^\infty z \, m(\dd z)
    < \infty ,
  \end{aligned}
 \end{equation}
 where we used that, by properties (ii) and (iii) of the function \ $\phi_k$, \ we have
 \ $\phi_k'(u) \in [0, 1]$ \ for all \ $u \in \RR$, \ and hence, by mean value theorem,
 \begin{equation}\label{(ii)}
  0 \leq \phi_k(y) - \phi_k(y - z) \leq z ,
  \qquad y \in \RR , \quad z \in \RR_+ , \quad k \in \NN .
 \end{equation}

One can check that the process \ $\left(I_{j,k,2}(t) + I_{j,k,4,1}(t)\right)_{t\in\RR_+}$
 \ is a martingale.
Indeed, by properties (ii) and (iii) of the function \ $\phi_k$ \ and the definition of
 \ $\tau_j$,
 \begin{align*}
  \EE\left( \int_0^{t\land\tau_j}
             \left(\phi_k'(Z_{s})
                   \bigl(\sqrt{Y'_s} - \sqrt{Y_s}\bigr)\right)^2
             \, \dd s  \right)
   \leq 2 \EE\left( \int_0^{t\land\tau_j} (Y'_s + Y_s) \, \dd s \right)
   \leq 4 j t < \infty ,
 \end{align*}
 hence, by Ikeda and Watanabe \cite[page 55]{IkeWat},
 \ $\left(I_{j,k,2}(t)\right)_{t\in\RR_+}$ \ is a martingale.
The inequality \eqref{Fp1} yields that the function
 \ $\RR_+ \times \RR_{++} \times \Omega \ni (s, z, \omega)
    \mapsto \phi_k(Z_{s-}(\omega) - z) - \phi_k(Z_{s-}(\omega)) \bone_{(-\infty,\tau_j)}(s)$
 \ belongs to the class \ $\bF_p^1$, \ and hence \ $\left(I_{j,k,4,1}(t)\right)_{t\in\RR_+}$
 \ is a martingale by Ikeda and Watanabe \cite[page 62]{IkeWat}.

Moreover,
 \begin{align*}
  I_{j,k,1}(t) \leq |b| \int_0^{t\land\tau_j} Z_s^+ \, \dd s .
 \end{align*}
By (iv),
 \[
   I_{j,k,3}(t) \leq (t\land\tau_j) \frac{\sigma^2}{k} \leq \frac{\sigma^2 t}{k} .
 \]
By \eqref{(ii)}, we obtain
 \[
   I_{j,k,4,2}(t) \leq 0 .
 \]
Summarizing, we have
 \begin{align}\label{help_comparison}
  \phi_k(Z_{t\land\tau_j})
  &\leq \phi_k(Z_0)
        + |b| \int_0^{t\land\tau_j} Z_s^+ \, \dd s
        + \frac{\sigma^2 t}{k}
        + I_{j,k,2}(t) + I_{j,k,4,1}(t) , \qquad t \in \RR_+ .
 \end{align}
By (iii), we obtain \ $\PP(\phi_k(Z_0) = 0) = 1$.
\ By (i), the non-negativeness of \ $\phi_k$ \ and monotone convergence
 theorem yield \ $\EE(\phi_k(Z_{t\land\tau_j})) \to \EE(Z_{t\land\tau_j}^+)$ \ as
 \ $k \to \infty$ \ for all \ $t \in \RR_+$ \ and \ $j \in \NN$.
\ We have
 \ $\int_0^{t\land\tau_j} Z_s^+ \, \dd s
    \leq \int_0^t Z_{s\land\tau_j}^+ \, \dd s$,
 \ hence taking the expectation of \eqref{help_comparison} and letting
 \ $k \to \infty$, \ we obtain
 \[
   \EE\bigl(Z_{t\land\tau_j}^+\bigr)
   \leq |b| \int_0^t \EE\bigl(Z_{s\land\tau_j}^+\bigr) \, \dd s .
 \]
By Gronwall's inequality, we conclude
 \[
   \EE\bigl(Z_{t\land\tau_j}^+\bigr) = 0
 \]
 for all \ $t \in \RR_+$ \ and \ $j \in \NN$.
\ Hence \ $\PP(Y'_{t\land\tau_j} \leq Y_{t\land\tau_j} )=1$ \ for all
 \ $t \in \RR_+$ \ and \ $j \in \NN$, \ and then
 \ $\PP(\text{$Y'_{t\land\tau_j} \leq Y_{t\land\tau_j}$ for all $j \in \NN$}) = 1$
 \ for all \ $t \in \RR_+$.
\ Since \ $Y$ \ and \ $Y'$ \ have c\`adl\`ag trajectories, these
 trajectories are bounded almost surely on \ $[0, T]$ \ for all
 \ $T \in \RR_+$, \ hence \ $\tau_j \as \infty$ \ as \ $j \to \infty$.
\ This yields \ $\PP(Y'_t \leq Y_t )=1$ \ for all \ $t \in \RR_+$.
\ Since the set of non-negative rational numbers \ $\QQ_+$ \ is countable, we
 obtain \ $\PP(\text{$Y'_t \leq Y_t$ for all $t \in \QQ_+$}) = 1$.
\ Using again that \ $Y$ \ and \ $Y'$ \ have c\`adl\`ag trajectories
 almost surely, we get
 \ $\PP(\text{$Y'_t \leq Y_t$ for all $t \in \RR_+$}) = 1$.
\proofend

\section{Likelihood-ratio process}
\label{App_LR}

Based on Jacod and Shiryaev \cite{JSh}, see also Jacod and M\'emin \cite{JM},
 S{\o}rensen \cite{SorM} and Luschgy \cite{Lus}, we recall certain sufficient
 conditions for the absolute continuity of probability measures induced by
 semimartingales together with a representation of the corresponding Radon--Nikodym
 derivative (likelihood-ratio process).
This appendix (together with proofs) already appears in Barczy et al. \cite{BarBenKebPap}, we
 decided to present it here as well for better readability and being self-contained.

Let \ $D(\RR_+, \RR^d)$ \ denote the space of
 \ $\RR^d$-valued c\`adl\`ag functions defined on \ $\RR_+$.
\ Let \ $(\eta_t)_{t\in\RR_+}$ \ denote the canonical process
 \ $\eta_t(\omega) := \omega(t)$, \ $\omega \in D(\RR_+, \RR^d)$,
 \ $t \in \RR_+$.
\ Put \ $\cF_t^\eta := \sigma(\eta_s, s \in [0, t])$, \ $t \in \RR_+$,
 \ and
 \[
   \cD_t(\RR_+, \RR^d)
   := \bigcap_{\vare \in \RR_{++}} \cF_{t+\vare}^\eta , \quad
   t \in \RR_+ , \qquad
   \cD(\RR_+, \RR^d)
   := \sigma\Biggl(\,\bigcup_{t \in \RR_+} \cF_t^\eta\Biggr) .
 \]
Let \ $\Psi \subset \RR^k$ \ be an arbitrary non-empty set, and let \ $\PP_\bpsi$,
 \ $\bpsi \in \Psi$, \ are probability measures on the canonical space
 \ $(D(\RR_+, \RR^d), \cD(\RR_+, \RR^d))$.
\ Suppose that for each \ $\bpsi \in \Psi$, \ under \ $\PP_\bpsi$, \ the canonical
 process \ $(\eta_t)_{t\in\RR_+}$ \ is a semimartingale with semimartingale
 characteristics \ $(B^{(\bpsi)}, C, \nu^{(\bpsi)})$ \ associated with a
 fixed Borel measurable truncation function \ $h:\RR^d \to \RR^d$, \ see Jacod and
 Shiryaev \cite[Definition II.2.6 and Remark II.2.8]{JSh}.
Namely, \ $C_t := \langle (\eta^\cont)^{(\bpsi)} \rangle_t$,
 \ $t \in \RR_+$, \ where
 \ $(\langle (\eta^\cont)^{(\bpsi)} \rangle_t)_{t\in\RR_+}$ \ denotes
 the (predictable) quadratic variation process (with values in \ $\RR^{d\times d}$)
 \ of the continuous martingale part \ $(\eta^\cont)^{(\bpsi)}$ \ of \ $\eta$ \ under
 \ $\PP_\bpsi$, \ $\nu^{(\bpsi)}$ \ is the compensator of the integer-valued random measure
 \ $\mu^\eta$ \ on \ $\RR_+ \times \RR^d$ \ associated with the jumps of
 \ $\eta$ \ under \ $\PP_\bpsi$ \ given by
 \[
   \mu^\eta(\omega, \dd t, \dd\bx)
   := \sum_{s\in\RR_+}
       \bbone_{\{\Delta\eta_s(\omega)\ne\bzero\}}
       \vare_{(s,\Delta\eta_s(\omega))}(\dd t, \dd\bx) , \qquad
   \omega \in D(\RR_+, \RR^d) ,
 \]
 where \ $\vare_{(t,\bx)}$ \ denotes the Dirac measure at the point
 \ $(t,\bx) \in \RR_+ \times \RR^d$, \ and
 \ $\Delta\eta_t := \eta_t - \eta_{t-}$, \ $t \in \RR_{++}$,
 \ $\Delta\eta_0 := \bzero$, \ and \ $B^{(\bpsi)}$ \ is the predictable process
 (with values in \ $\RR^d$ \ having finite variation over each finite interval
 \ $[0, t]$, \ $t \in \RR_+$) \ appearing in the canonical decomposition
 \[
   \teta_t = \eta_0 + M^{(\bpsi)}_t + B^{(\bpsi)}_t , \qquad
   t \in \RR_+ ,
 \]
 of the special semimartingale \ $(\teta_t)_{t\in\RR_+}$ \ under \ $\PP_\bpsi$
 \ given by
 \[
   \teta_t := \eta_t - \sum_{s\in[0,t]} (\eta_s - h(\Delta\eta_s)) , \qquad
   t \in \RR_+ ,
 \]
 where \ $(M^{(\bpsi)}_t)_{t\in\RR_+}$ \ is a local martingale with
 \ $M^{(\bpsi)}_0 = \bzero$.
\ We call the attention that, by our assumption, the process
 \ $C= \langle (\eta^\cont)^{(\bpsi)} \rangle$ \ does not depend on \ $\bpsi$,
 \ although \ $(\eta^\cont)^{(\bpsi)}$ \ might depend on \ $\bpsi$.
\ In addition, assume that \ $\PP_\bpsi(\nu^{(\bpsi)}(\{t\} \times \RR^d) = 0)=1$
 \ for every \ $\bpsi \in \Psi$, \ $t \in \RR_+$, \ and
 \ $\PP_\bpsi(\eta_0 = \bx_0) = 1$ \ with some \ $\bx_0 \in \RR^d$
 \ for every \ $\bpsi \in \Psi$.
\ Note that we have the semimartingale representation
 \begin{align*}
  \eta_t &= \bx_0 + B^{(\bpsi)}_t + (\eta^\cont)^{(\bpsi)}_t
            + \int_0^t \int_{\RR^d} h(\bx) \,
               (\mu^\eta - \nu^{(\bpsi)})(\dd s, \dd\bx) \\
         &\quad
            + \int_0^t \int_{\RR^d} (\bx - h(\bx)) \, \mu^\eta(\dd s, \dd\bx) ,
  \qquad t \in \RR_+ ,
 \end{align*}
 of \ $\eta$ \ under \ $\PP_\bpsi$, \ see Jacod and Shiryaev
 \cite[Theorem II.2.34]{JSh}.
Moreover, for each \ $\bpsi \in \Psi$, \ let us choose a nondecreasing,
 continuous, adapted process \ $(F_t^{(\bpsi)})_{t\in\RR_+}$ \ with
 \ $F_0^{(\bpsi)} = 0$ \ and a predictable process
 \ $(c_t^{(\bpsi)})_{t\in\RR_+}$ \ with values in the set of all
 symmetric positive semidefinite \ $d \times d$ \ matrices such that
 \[
   C_t = \int_0^t c_s^{(\bpsi)} \, \dd F_s^{(\bpsi)}
 \]
 $\PP_\bpsi$-almost sure for every \ $t \in \RR_+$.
\ Due to the assumption \ $\PP_\bpsi(\nu^{(\bpsi)}(\{t\} \times \RR^d) = 0) = 1$
 \ for every \ $\bpsi \in \Psi$, \ $t \in \RR_+$, \ such choices of
 \ $(F_t^{(\bpsi)})_{t\in\RR_+}$ \ and
 \ $(c_t^{(\bpsi)})_{t\in\RR_+}$ \ are possible, see Jacod and
 Shiryaev \cite[Proposition II.2.9 and Corollary II.1.19]{JSh}.
Let \ $\cP$ \ denote the predictable $\sigma$-algebra on
 \ $D(\RR_+, \RR^d) \times \RR_+$.
\ Assume also that for every \ $\bpsi, \btpsi \in \Psi$, \ there exist a
 \ $\cP \otimes \cB(\RR^d)$-measurable function
 \ $V^{(\btpsi,\bpsi)}
    : D(\RR_+, \RR^d) \times \RR_+ \times \RR^d \to \RR_{++}$
 \ and a predictable \ $\RR^d$-valued process \ $\beta^{(\btpsi,\bpsi)}$
 \ satisfying
 \begin{gather}
  \nu^{(\bpsi)}(\dd t, \dd\bx)
  = V^{(\btpsi,\bpsi)}(t, \bx) \nu^{(\btpsi)}(\dd t, \dd\bx) , \label{GIR1}\\
  \int_0^t \int_{\RR^d}
   \Bigl(\sqrt{V^{(\btpsi,\bpsi)}(s, \bx)} - 1\Bigr)^2 \,
   \nu^{(\btpsi)}(\dd s, \dd\bx)
  < \infty , \label{GIR2} \\
  B^{(\bpsi)}_t
  = B^{(\btpsi)}_t
    + \int_0^t
       c_s^{(\bpsi)} \beta^{(\btpsi,\bpsi)}_s \, \dd F_s^{(\bpsi)}
    + \int_0^t \int_{\RR^d}
       (V^{(\btpsi,\bpsi)}(s, \bx) - 1) h(\bx) \,
       \nu^{(\btpsi)}(\dd s, \dd\bx) , \label{GIR3} \\
  \int_0^t
   (\beta^{(\btpsi,\bpsi)}_s)^\top c_s^{(\bpsi)} \beta^{(\btpsi,\bpsi)}_s \, \dd F_s^{(\bpsi)}
  < \infty , \label{GIR4}
 \end{gather}
 $\PP_\bpsi$-almost sure for every \ $t \in \RR_+$.
\ Further, assume that for each \ $\bpsi \in \Psi$, \ local uniqueness holds for the
 martingale problem on the canonical space corresponding to the triplet
 \ $(B^{(\bpsi)}, C, \nu^{(\bpsi)})$ \ with the given initial value \ $\bx_0$
 \ with \ $\PP_\bpsi$ \ as its unique solution.
Then for each \ $T \in \RR_+$, \ $\PP_{\bpsi,T}$ \ is absolutely continuous
 with respect to \ $\PP_{\btpsi,T}$,
 \ where \ $\PP_{\bpsi,T} := \PP_{\bpsi\!\!}|_{\cD_T(\RR_+, \RR^d)}$ \ denotes the
 restriction of \ $\PP_\bpsi$ \ to \ $\cD_T(\RR_+, \RR^d)$ \ (similarly for
 \ $\PP_{\btpsi,T}$), \ and, under \ $\PP_{\btpsi,T}$, \ the corresponding
 likelihood-ratio process takes the form
 \begin{align}\label{RN_general}
  \begin{split}
  \log \frac{\dd \PP_{\bpsi,T}}{\dd \PP_{\btpsi,T}}(\eta)
  &= \int_0^T (\beta^{(\btpsi,\bpsi)}_s)^\top \, \dd(\eta^\cont)^{(\btpsi)}_s
     - \frac{1}{2}
       \int_0^T
       (\beta^{(\btpsi,\bpsi)}_s)^\top c_s^{(\bpsi)} \beta^{(\btpsi,\bpsi)}_s \, \dd F_s^{(\bpsi)} \\
  &\quad
    + \int_0^T \int_{\RR^d}
       (V^{(\btpsi,\bpsi)}(s, \bx) - 1) \,
       (\mu^\eta - \nu^{(\btpsi)})(\dd s, \dd\bx)\\
  &\quad
    + \int_0^T \int_{\RR^d}
       (\log(V^{(\btpsi,\bpsi)}(s, \bx))
        - V^{(\btpsi,\bpsi)}(s, \bx) + 1) \,
       \mu^\eta(\dd s, \dd\bx)
  \end{split}
 \end{align}
 for all \ $T \in \RR_{++}$, \ see Jacod and Shiryaev \cite[Theorem III.5.34]{JSh}.
A detailed proof of \eqref{RN_general} using Jacod and Shiryaev \cite{JSh} can be
 found in Barczy et al.\ \cite[Appendix A]{BarBenKebPap}.

\section{Limit theorems for continuous local martingales}
\label{App_clm}

In what follows we recall some limit theorems for continuous local
 martingales.
We use these limit theorems for studying the asymptotic
 behaviour of the MLE of \ $b$.
\ First we recall a strong law of large numbers for continuous local
 martingales.

\begin{Thm}{\bf (Liptser and Shiryaev \cite[Lemma 17.4]{LipShiII})}
\label{DDS_stoch_int}
Let \ $\bigl( \Omega, \cF, (\cF_t)_{t\in\RR_+}, \PP \bigr)$ \ be a filtered
 probability space satisfying the usual conditions.
Let \ $(M_t)_{t\in\RR_+}$ \ be a square-integrable continuous local martingale with respect to the
 filtration \ $(\cF_t)_{t\in\RR_+}$ \ such that \ $\PP(M_0 = 0) = 1$.
\ Let \ $(\xi_t)_{t\in\RR_+}$ \ be a progressively measurable process such that
 \[
   \PP\left( \int_0^t \xi_u^2 \, \dd \langle M \rangle_u < \infty \right) = 1 ,
   \qquad t \in \RR_+ ,
 \]
 and
 \begin{align}\label{SEGED_STRONG_CONSISTENCY2}
  \int_0^t \xi_u^2 \, \dd \langle M \rangle_u \as \infty \qquad
  \text{as \ $t \to \infty$,}
 \end{align}
 where \ $(\langle M \rangle_t)_{t\in\RR_+}$ \ denotes the quadratic variation process of
 \ $M$.
\ Then
 \begin{align}\label{SEGED_STOCH_INT_SLLN}
  \frac{\int_0^t \xi_u \, \dd M_u}
       {\int_0^t \xi_u^2 \, \dd \langle M \rangle_u} \as 0 \qquad
  \text{as \ $t \to \infty$.}
 \end{align}
If \ $(M_t)_{t\in\RR_+}$ \ is a standard Wiener process, the progressive measurability of
 \ $(\xi_t)_{t\in\RR_+}$ \ can be relaxed to measurability and adaptedness to the filtration
 \ $(\cF_t)_{t\in\RR_+}$.
\end{Thm}

The next theorem is about the asymptotic behaviour of continuous multivariate local
 martingales, see van Zanten \cite[Theorem 4.1]{Zan}.

\begin{Thm}{\bf (van Zanten \cite[Theorem 4.1]{Zan})}\label{THM_Zanten}
Let \ $\bigl( \Omega, \cF, (\cF_t)_{t\in\RR_+}, \PP \bigr)$ \ be a filtered
 probability space satisfying the usual conditions.
Let \ $(\bM_t)_{t\in\RR_+}$ \ be a $d$-dimensional square-integrable continuous local martingale
 with respect to the filtration \ $(\cF_t)_{t\in\RR_+}$ \ such that
 \ $\PP(\bM_0 = \bzero) = 1$.
\ Suppose that there exists a function \ $\bQ : \RR_+ \to \RR^{d \times d}$
 \ such that \ $\bQ(t)$ \ is an invertible (non-random) matrix for all
 \ $t \in \RR_+$, \ $\lim_{t\to\infty} \|\bQ(t)\| = 0$ \ and
 \[
   \bQ(t) \langle \bM \rangle_t \, \bQ(t)^\top \stoch \bfeta \bfeta^\top
   \qquad \text{as \ $t \to \infty$,}
 \]
 where \ $\bfeta$ \ is a \ $d \times d$ random matrix.
Then, for each $\RR^k$-valued random vector \ $\bv$ \ defined on
 \ $(\Omega, \cF, \PP)$, \ we have
 \[
   (\bQ(t) \bM_t, \bv) \distr (\bfeta \bZ, \bv) \qquad
   \text{as \ $t \to \infty$,}
 \]
 where \ $\bZ$ \ is a \ $d$-dimensional standard normally distributed random
 vector independent of \ $(\bfeta, \bv)$.
\end{Thm}

We note that Theorem \ref{THM_Zanten} remains true if the function \ $\bQ$
 \ is defined only on an interval \ $[t_0, \infty)$ \ with some
 \ $t_0 \in \RR_{++}$.

\section*{Acknowledgements}
We would like to thank the referee for his/her comments that helped us to improve the paper.


\begin{thebibliography}{99}

\bibitem{Alf}
\textsc{Alfonsi, A.}, (2005).
On the discretization schemes for the CIR (and Bessel squared) processes.
\textit{Monte Carlo Methods and Applications}
\textbf{11(4)} 355--384.

\bibitem{BarBenKebPap}
\textsc{Barczy, M., Ben Alaya, M., Kebaier, A.} and \textsc{Pap, G.} (2016).
Asymptotic behavior of maximum likelihood estimators for a jump-type Heston model.
Available on ArXiv: \url{http://arxiv.org/abs/1509.08869}

\bibitem{BarDorLiPap}
\textsc{Barczy, M.}, \textsc{D\"oring, L.}, \textsc{Li, Z.} and
\textsc{Pap, G.} (2013).
On parameter estimation for critical affine processes.
\textit{Electronic Journal of Statistics}
\textbf{7} 647--696.

\bibitem{BarLiPap1}
\textsc{Barczy, M.}, \textsc{Li, Z.} and \textsc{Pap, G.} (2015).
Yamada--Watanabe results for stochastic differential equations with jumps.
\textit{International Journal of Stochastic Analysis}
\textbf{Volume 2015} Article ID 460472, 23 pages.

\bibitem{BarLiPap2}
\textsc{Barczy, M., Li, Z.} and \textsc{Pap, G.} (2015).
Stochastic differential equation with jumps for multi-type continuous state and continuous
 time branching processes with immigration.
\textit{ALEA. Latin American Journal of Probability and Mathematical Statistics}.
\textbf{12(1)} 129--169.

\bibitem{BenKeb1}
\textsc{Ben Alaya, M.} and \textsc{Kebaier, A.} (2012).
Parameter estimation for the square root diffusions: ergodic and nonergodic cases.
\textit{Stochastic Models}
\textbf{28(4)} 609--634.

\bibitem{BenKeb2}
\textsc{Ben Alaya, M.} and \textsc{Kebaier, A.} (2013).
Asymptotic behavior of the maximum likelihood estimator for ergodic and nonergodic
 square-root diffusions.
\textit{Stochastic Analysis and Applications}
\textbf{31(4)} 552--573.

\bibitem{Bil}
\textsc{Billingsley, P.} (1999).
\textit{Convergence of probability measures}, 2nd ed.
John Wiley \& Sons, Inc., New York.

\bibitem{Bha}
\textsc{Bhattacharya, R. N.} (1982).
On the functional central limit theorem and the law of the iterated logarithm for Markov
 processes.
\textit{Zeitschrift f\"ur Wahrscheinlichkeitstheorie und Verwandte Gebiete}
\textbf{60} 185--201.

\bibitem{CoxIngRos}
\textsc{Cox, J. C.}, \textsc{Ingersoll, J. E.} and \textsc{Ross, S. A.} (1985).
A theory of the term structure of interest rates.
\textit{Econometrica}
\textbf{53(2)} 385--407.

\bibitem{DawLi}
\textsc{Dawson, D. A.} and \textsc{Li, Z.} (2006).
Skew convolution semigroups and affine Markov processes.
\textit{The Annals of Probability}
\textbf{34(3)} 1103--1142.

\bibitem{Dud}
\textsc{Dudley, R. M.} (1989).
\textit{Real Analysis and Probability}.
Wadsworth \& Brooks/Cole Advanced Books \& Software, Pacific Grove, California.

\bibitem{DufGar}
\textsc{Duffie, D.} and \textsc{G\^arleanu, N.} (2001).
Risk and valuation of collateralized debt obligations.
\textit{Financial Analysts Journal}
\textbf{57(1)} 41--59.

\bibitem{DufFilSch}
\textsc{Duffie, D.}, \textsc{Filipovi\'{c}, D.} and \textsc{Schachermayer, W.} (2003).
Affine processes and applications in finance.
\textit{Annals of Applied Probability}
\textbf{13} 984--1053.

\bibitem{Fel}
\textsc{Feller, W.} (1951).
Two singular diffusion problems.
\textit{Annals of Mathematics}
\textbf{54(1)} 173--182.

\bibitem{Fil}
\textsc{Filipovi\'c, D.} (2001).
A general characterization of one factor affine term structure models.
\textit{Finance and Stochastics}
\textbf{5(3)} 389--412.

\bibitem{FilMaySch}
\textsc{Filipovi\'{c}, D., Mayerhofer, E.} and \textsc{Schneider, P.} (2013).
Density approximations for multivariate affine jump-diffusion processes.
\textit{Journal of Econometrics}
\textbf{176(2)} 93--111.

\bibitem{FuLi}
\textsc{Fu, Z.} and \textsc{Li, Z.} (2010).
Stochastic equations of non-negative processes with jumps.
\textit{Stochastic Processes and their Applications}
\textbf{120} 306--330.

\bibitem{HuaMaZhu}
\textsc{Huang, J., Ma, C.} and \textsc{Zhu, C.} (2011).
Estimation for discretely observed continuous state branching processes with immigration.
\textit{Statistics and Probability Letters}
\textbf{81} 1104--1111.

\bibitem{IkeWat}
\textsc{Ikeda, N.} and \textsc{Watanabe, S.} (1989).
\textit{Stochastic Differential Equations and Diffusion Processes}, 2nd ed.
North-Holland/Kodansha, Amsterdam/Tokyo.

\bibitem{JM}
\textsc{Jacod, J.} and \textsc{M\'emin, J.} (1976).
Caract\'eristiques locales et conditions de continuit\'e absolue pour les semi-martingales.
\textit{Zeitschrift f\"ur Wahrscheinlichkeitstheorie und Verwandte Gebiete}, \textbf{35}
 1--37.

\bibitem{JacPro}
\textsc{Jacod, J.} and \textsc{Protter, P.} (2012).
\textit{Discretization of processes},
Springer, Heidelberg.

\bibitem{JSh}
\textsc{Jacod, J.} and \textsc{Shiryaev, A. N.} (2003).
\textit{Limit Theorems for Stochastic Processes}, 2nd ed.
Springer-Verlag, Berlin.

\bibitem{JiaMaSco}
\textsc{Jiao, Y., Ma, C.} and \textsc{Scotti, S.} (2016).
Alpha-CIR model with branching processes in sovereign interest rate modelling.
Available on ArXiv: \url{http://arxiv.org/abs/1602.05541}

\bibitem{JinRudTra2}
\textsc{Jin, P., R\"udiger, B.} and \textsc{Trabelsi, C.} (2016).
Exponential ergodicity of the jump-diffusion CIR process.
In: Stochastics of environmental and financial economics---Centre of Advanced Study,
Oslo, Norway, 2014--2015, pages 285--300,
Springer Proceedings in Mathematics \& Statistics, 138, Springer, Cham.

\bibitem{KR2}
\textsc{Keller-Ressel, M.} (2008).
 Affine Processes - Theory and Applications in Finance.
\textit{PhD Thesis, Vienna University of Technology,} pages 110.

\bibitem{KR}
\textsc{Keller-Ressel, M.} (2011).
Moment explosions and long-term behavior of affine stochastic volatility models.
\textit{Mathematical Finance}
\textbf{21(1)} 73--98.

\bibitem{KRMij}
\textsc{Keller-Ressel, M.} and \textsc{Mijatovi\'c, A.} (2012).
On the limit distributions of continuous-state branching processes with immigration.
\textit{Stochastic Processes and their Applications}
\textbf{122} 2329--2345.

\bibitem{KRSte}
\textsc{Keller-Ressel, M.} and \textsc{Steiner, T.} (2008).
Yield curve shapes and the asymptotic short rate distribution in affine one-factor models.
\textit{Finance and Stochastics}
\textbf{12(2)} 149--172.

\bibitem{KucSor}
\textsc{K\"uchler, U.} and \textsc{S{\o}rensen, M.} (1997).
\textit{Exponential families of stochastic processes,}
Springer-Verlag, New York.

\bibitem{Kyp}
\textsc{Kyprianou, A. E.} (2014).
\textit{Fluctuations of L\'evy Processes with Applications,} 2nd ed.
 Springer-Verlag, Berlin Heidelberg.

\bibitem{LamLap}
\textsc{Lamberton, D.} and \textsc{Lapeyre, B.} (1996).
\textit{Introduction to Stochastic Calculus Applied to Finance}.
Chapman \& Hall/CRC.

\bibitem{Li}
\textsc{Li, Z.} (2011).
\textit{Measure-Valued Branching Markov Processes}.
Springer-Verlag, Heidelberg.

\bibitem{Li2}
\textsc{Li, Z.} (2012).
\textit{Continuous-state branching processes}.
Available on ArXiv: \url{http://arxiv.org/abs/1202.3223}

\bibitem{LiMa}
\textsc{Li, Z.} and \textsc{Ma, C.} (2013).
 Asymptotic properties of estimators in a stable Cox-Ingersoll-Ross model.
\textit{Stochastic Processes and their Applications} \textbf{125(8)} 3196--3233.

\bibitem{LipShiII}
\textsc{Liptser, R. S.} and \textsc{Shiryaev, A. N.} (2001).
\textit{Statistics of Random Processes II. Applications}, 2nd edition.
Springer-Verlag, Berlin, Heidelberg.

\bibitem{Lus2}
\textsc{Luschgy, H.} (1992).
Local asymptotic mixed normality for semimartingale experiments.
\textit{Probability Theory and Related Fields}
\textbf{92} 151--176.


\bibitem{Lus}
\textsc{Luschgy, H.} (1994).
Asymptotic inference for semimartingale models with singular parameter points.
\textit{Journal of Statistical Planning and Inference}
\textbf{39} 155--186.

\bibitem{Ma}
\textsc{Ma, R.} (2013).
Stochastic equations for two-type continuous-state branching processes with immigration.
\textit{Acta Mathematica Sinica, English Series}
\textbf{29(2)} 287--294.

\bibitem{Mai}
\textsc{Mai, H.} (2012).
Drift estimation for jump diffusions: time-continuous and high-frequency observations.
Ph.D. Dissertation. \textit{Humboldt-Universit\"{a}t zu Berlin.}

\bibitem{Ove}
\textsc{Overbeck, L.} (1998).
Estimation for continuous branching processes.
\textit{Scandinavian Journal of Statistics}
\textbf{25(1)} 111--126.

\bibitem{OveRyd}
\textsc{Overbeck, L.} and \textsc{Ryd\'en, T.} (1997).
Estimation in the Cox-Ingersoll-Ross model.
\textit{Econometric Theory}
\textbf{13(3)} 430--461.

\bibitem{Pin}
\textsc{Pinsky, M. A.} (1972).
Limit theorems for continuous state branching processes with immigration.
\textit{Bulletin of the American Mathematical Society}
\textbf{78(2)} 242--244.

\bibitem{Sat}
\textsc{Sato, K.-I.} (1999).
\textit{L\'evy Processes and Infinitely Divisible Distributions.}
 Cambridge University Press, Cambridge.

\bibitem{SorM}
\textsc{S{\o}rensen, M.} (1991).
Likelihood methods for diffusions with jumps.
In: N. U. Prabhu and I.V. Basawa, Eds.,
\textit{Statistical Inference in Stochastic Processes,}
Marcel Dekker, New York, 67--105.

\bibitem{Zan}
\textsc{van Zanten, H.} (2000).
A multivariate central limit theorem for continuous local martingales.
\textit{Statistics \& Probability Letters}
\textbf{50(3)} 229--235.

\end{thebibliography}
\end{document}